\documentclass[11pt]{article}
\usepackage{amscd}
\usepackage{amsmath}
\usepackage{latexsym}
\usepackage{amsfonts}
\usepackage{amssymb}
\usepackage{amsthm}
\usepackage{graphicx}

\newtheorem{theorem}{Theorem}[section]
\newtheorem{lemma}{Lemma}[section]

\newtheorem{proposition}{Proposition}[section]
\newtheorem{definition}{Definition}[section]

\newtheorem{corollary}{Corollary}[section]
\newtheorem{remark}{Remark}[section]

\makeatletter
\newcommand{\Extend}[5]{\ext@arrow0099{\arrowfill@#1#2#3}{#4}{#5}}
\makeatother

\begin{document}

 \title{ Global well-posedness, scattering and blow-up for the energy-critical, focusing Hartree equation in the radial case}
 \author{{Changxing Miao,\ \ Guixiang Xu,\ \ and \ Lifeng Zhao }\\
         {\small Institute of Applied Physics and Computational Mathematics}\\
         {\small P. O. Box 8009,\ Beijing,\ China,\ 100088}\\
         {\small (miao\_changxing@iapcm.ac.cn, \ xu\_guixiang@iapcm.ac.cn, zhao\_lifeng@iapcm.ac.cn ) }\\
         \date{}
        }
\maketitle

\begin{abstract}
We establish global existence, scattering for radial solutions to
the  energy-critical  focusing Hartree equation with energy and
$\dot{H}^1$ norm less than those of the ground state in
$\mathbb{R}\times \mathbb{R}^d$, $d\geq 5$.
\end{abstract}

 \begin{center}
 \begin{minipage}{120mm}
 \end{minipage}
 \end{center}

%{\center \tableofcontents}

\section{Introduction}
 \setcounter{section}{1}\setcounter{equation}{0}

We consider the following initial value problem
\begin{equation} \label{equ}
\left\{ \aligned
    iu_t +  \Delta u  & = f(u), \quad  \text{in}\  \mathbb{R}^d \times \mathbb{R}, \quad d\geq 5,\\
     u(0)&=u_0(x), \ \ \text{in} \ \mathbb{R}^d,
\endaligned
\right.
\end{equation}
where $u(t,x)$ is a complex-valued function in spacetime $\mathbb{R}
\times \mathbb{R}^d$ and $\Delta$ is the Laplacian in
$\mathbb{R}^{d}$, $f(u)=-\big( |x|^{-4}* |u|^2 \big) u$. It is
introduced as a classical model in \cite{website}. In practice, we
use the integral formulation of $(\ref{equ})$
\begin{equation}\label{intequa1}
u(t)=U(t)u_0(x) -i \int^{t}_{0} U(t-s)f(u(s))ds,
\end{equation}
where $U(t)=e^{it\Delta}$.

%This equation has the Hamiltonian
%\begin{equation}
%\aligned E(u(t))=\frac12 & \big\|\nabla
%u(t)\big\|^2_{L^2}-\frac{1}{4} \iint \frac{1}{|x-y|^{4}}
%|u(t,x)|^2 |u(t,y)|^2\ dxdy.
%\endaligned
%\end{equation}
%Since $(\ref{energy})$ is preserved by the flow corresponding to
%$(\ref{equ1})$ we shall refer to it as the energy and often write
%$E(u)$ for $E(u(t))$.

%A second conserved quantity we will occasionally rely on is the mass
%$\big\|u(t)\big\|_{L^2_x(\mathbb{R}^n)}$. However, since the
%equation is $L^2$-supercritical with respect to the scaling (see
%$(\ref{scaling})$), we do not have bounds on the mass that are
%uniform across frequencies. Indeed, the low frequencies may
%simultaneously have small energy and large mass.

We are primarily interested in $(\ref{equ})$ since it is critical
with respect to the energy norm. That is, the scaling $u\mapsto
u_{\lambda}$ where
\begin{equation}\label{scaling}
u_{\lambda}(t,x)=\lambda^{\frac{d-2}{2}}u(\lambda^2 t, \lambda x), \
\lambda>0
\end{equation}
maps a solution to $(\ref{equ})$ to another solution to
$(\ref{equ})$, and $u$ and $u_{\lambda}$ have the same energy (\ref{energy}).

It is known that if the initial data $u_0(x)$ has finite energy,
then $(\ref{equ})$ is locally well-posed (see, for instance
\cite{MiXZ06}). That is, there exists a unique local-in-time
solution that lies in $C^0_t\dot{H}^1_x \cap
L^6_tL^{\frac{6d}{3d-8}}_x$ and the map from the initial data to the
solution is locally Lipschitz in these norms. If the energy is
small, it is known that the solution exists globally in time and
scattering occurs; That is, there exist solutions $u_{\pm}$ of the
free Schr\"{o}dinger equation $(i\partial_t + \Delta)u_{\pm}=0$ such
that
\begin{equation*}
\aligned
 \big\|u(t)-u_{\pm}(t)\big\|_{\dot{H}^1_x} \rightarrow 0 \quad \text{as}\ t\rightarrow \pm \infty.
 \endaligned
 \end{equation*}
However, for initial data with large
energy, the local well-posedness argument do not extend to give
global well-posedness, only with the conservation of the energy
$(\ref{energy})$, because the time of existence given by the local
theory depends on the profile of the data as well as on $\big\|u_0\big\|_{\dot{H}^1_x}$.

A large amount of work has been devoted to the theory of scattering
for the Hartree equation, see \cite{GiO93}-\cite{HaT87},
\cite{Mi97}-\cite{MiXZ07b}, \cite{Na99d} and \cite{NaO92}. In particular, global
well-posedness in $\dot{H}^1_x$ for the energy-critical, defocusing
Hartree equation in the case of large finite-energy initial data was obtained recently by us
\cite{MiXZ07a}, \cite{MiXZ07b}. In this paper, we continue this investigation and
establish scattering result for radial solutions to the energy-critical, focusing Hartree equation for data with
energy and $\dot{H}^1$ norm less than those of the gound state.

\begin{figure}[h]
\qquad \qquad \qquad \qquad \qquad  \includegraphics[width=0.5\textwidth, angle=0]{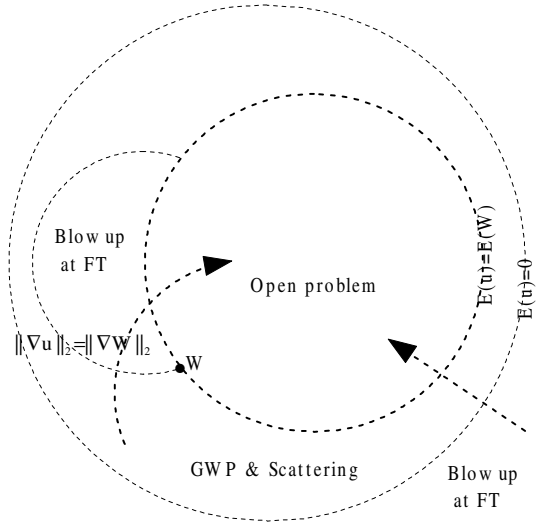}
\caption[]{A description of the solutions with radial data in the energy space, where ``FT" denotes finite time.}
\end{figure}

The main result of this paper is the following global well-posedness and blow up results
for $(\ref{equ})$ in the energy space (Figure 1).

\begin{theorem}\label{main}
Let $d\geq 5$, $u_0 \in \dot{H}^1(\mathbb{R}^d)$ be radial and let $u$ be the corresponding solution to $(\ref{equ})$ in $\dot{H}^1(\mathbb{R}^d)$ with maximal forward time
interval of existence $[0, T)$. Suppose $E(u_0) < E(W)$.
\begin{enumerate}
\item[$(1)$] If $\big\|\nabla u_0\big\|_{L^2}<\big\|\nabla W\big\|_{L^2}$, then $T=+\infty$ and $u$ scatters in $\dot{H}^1$.
\item[$(2)$] If $\big\|\nabla u_0\big\|_{L^2}>\big\|\nabla W\big\|_{L^2}$,
then $T<+\infty$, and thus, the solution blows up at finite time.
\end{enumerate}
\end{theorem}

Similar as in \cite{KeM06}, it is still open that scattering for the general data with energy and $\dot{H}^1$ norm less than those of the gound state.
But concerning the blow up result, we also have
\begin{theorem}\label{blowupresult2}
Let $d\geq 5$, $u_0 \in \dot{H}^1(\mathbb{R}^d)$ and let $u$ be the corresponding solution to $(\ref{equ})$ in $\dot{H}^1(\mathbb{R}^d)$ with maximal forward time
interval of existence $[0, T)$. Suppose $E(u_0) < E(W)$, $\big\|\nabla u_0\big\|_{L^2}>\big\|\nabla W\big\|_{L^2}$ and $|x|u_0 \in L^2$,
then $T<+\infty$, i.e., the solution blows up at finite time.
\end{theorem}

%The results here have analogs in previous works on nonlinear
%Schr\"{o}dinger equations. Global wellposedness and scattering in
%$\mathbb{R}^{1+3}$ were shown by I-team \cite{CKSTT07}, global
%wellposedness and scattering in $\mathbb{R}^{1+4}$ were shown by
%Ryckman-Visan \cite{RyV05}, and global wellposedness and scattering
%in $\mathbb{R}^{1+n}$ ($n\geq 5$) were shown by Visan \cite{Vi05}
%and \cite{Vi06}. The result of Ryckman and Visan corresponds to the
%$n=4$ limit of our result. Because the nonlinearity is the nonlocal
%term, we need the dyadic decomposition and a new Bernstein estimate
%with the stability theory to set up the frequency localization of
%the minimal energy blowup solution (definition in section 4) at each
%time. On the other hand, according to our needs, we prove that the
%$L^{ 2n/n-2}_x$-norm of the minimal energy blowup solution is
%bounded from below, instead that the potential energy is bounded
%from below.

Next, we introduce some notations. If $X, Y$ are nonnegative quantities, we use $X\lesssim Y $ or
$X=O(X)$ to denote the estimate $X\leq CY$ for some $C$ which may
depend on the critical energy $E_{crit}$ (see Section 4) but not on
any parameter such as $\eta$, and $X \approx Y$ to denote the
estimate $X\lesssim Y\lesssim X$. We use $X\ll Y$ to mean $X \leq c
Y$ for some small constant $c$ which is again allowed to depend on
$E_{crit}$.

We use $C\gg1$ to denote various large finite constants. and $0< c
\ll 1$ to denote various small constants.

The Fourier transform on $\mathbb{R}^d$ is defined by
\begin{equation*}
\aligned \widehat{f}(\xi):= \big( 2\pi
\big)^{-\frac{d}{2}}\int_{\mathbb{R}^d}e^{- ix\cdot \xi}f(x)dx ,
\endaligned
\end{equation*}
giving rise to the fractional differentiation operators
$|\nabla|^{s}$,  defined by
\begin{equation*}
\aligned \widehat{|\nabla|^sf}(\xi):=|\xi|^s\widehat{f}(\xi).
\endaligned
\end{equation*}
These define the homogeneous Sobolev norms
\begin{equation*}
\big\|f\big\|_{\dot{H}^s_x}:= \big\| |\nabla|^sf
\big\|_{L^2_x(\mathbb{R}^d)}.
\end{equation*}

Let $e^{it\Delta }$ be the free Schr\"{o}dinger propagator. In
physical space this is given by the formula
\begin{equation*}
\aligned e^{it\Delta }f(x)=\frac{1}{(4\pi it)^2}\int_{\mathbb{R}^d}
e^{\frac{i|x-y|^2}{4t}}f(y)dy,
\endaligned
\end{equation*}
while in frequency space one can write this as
\begin{equation*}
\widehat{e^{it\Delta } f}(\xi)=e^{- it|\xi|^2}\widehat{f}(\xi).
\end{equation*}

In particular, the propagator preserves the above Sobolev norms and
obeys the dispersive estimate
\begin{equation}\label{dispersivee}
\aligned \big\|e^{it\Delta} f\big\|_{L^{\infty}_x(\mathbb{R}^d)}
\lesssim |t|^{-\frac{d}{2}}\big\|f\big\|_{L^1_x(\mathbb{R}^d)}, \quad \forall \ t\not=0 .
\endaligned
\end{equation}

Let $d\geq 5$, a pair $(q, r)$ is $L^2$-admissible if
\begin{equation*}
\frac{2}{q} = d\Big(\frac{1}{2}-\frac{1}{r}\Big), \ \text{for} \ \
2\leq r \leq \frac{2d}{d-2}.
\end{equation*}

For a spacetime slab $I\times \mathbb{R}^d$, we define the {\it
Strichartz} norm $\dot{S}^0(I)$ by
\begin{equation*}
\big\|u\big\|_{\dot{S}^0(I)}:= \sup_{(q, r)\ \text{$L^2$-admissible}}
\big\|u\big\|_{L^q_tL^r_x(I\times \mathbb{R}^d)}.
\end{equation*}
and for some fixed number $0< \epsilon_0\ll 1$, define ${\cal Z}^1(I)$ by
\begin{equation*}
\big\|u\big\|_{{\cal Z}^1(I)}:= \sup_{(q, r) \in \wedge }\big\| u\big\|_{L^q_tL^r_x},
\end{equation*}
where
\begin{equation*}
\aligned
 \wedge =\Big\{(q,r); \frac2q=d(\frac12-\frac1r)-1, \frac{2d}{d-2}\leq r \leq \frac{2d}{d-4} -\epsilon_0 \Big\}.
\endaligned
\end{equation*}

When $d\geq 5$, the spaces $\big(\dot{S}^0(I),
\|\cdot\|_{\dot{S}^0(I)}\big)$ and $\big({\cal Z}^1(I),
\|\cdot\|_{{\cal Z}^1(I)}\big)$ are Banach spaces, respectively.

We will occasionally use subscripts to denote spatial derivatives
and will use the summation convention over repeated indices.

We work in the frame of \cite{KeM06}, \cite{KeM06-2} and \cite{KiTV07}. In Section $2$, we recall some
useful facts. In Section $3$, we obtain some variational estimates and blow up results (Part $(2)$ of Theorem \ref{main} and
Theorem \ref{blowupresult2}). Last using a concentration compactness argument, we obtain the scattering result (Part $(1)$ of Theorem \ref{main}) in Section $4$ and $5$.

%%%%%%%%%%%%%%%%%%%%%%%%%%%%%%%%%%%%%%%%%%%%%%%%%%%%%%%%%%%%%%%%%%%%%%%%%%%%%%%%%%%
\section{A review of the Cauchy problem}
 \setcounter{section}{2}
\setcounter{equation}{0} In this section, we will recall some basic facts about the Cauchy
problem
\begin{equation} \label{equ1}
\left\{ \aligned
    iu_t +  \Delta u  & = f(u),  & (x, t) \in \mathbb{R}^d \times \mathbb{R}, \quad d\geq 5,\\
     u(t_0) \in & \ \dot{H}^1(\mathbb{R}^d),
\endaligned
\right.
\end{equation}
where $f(u)=-\big( |x|^{-4}* |u|^2 \big) u$. It is the $\dot{H}^1$
critical, focusing Hartree equation.

Based on the above notations, we have the following {\it Strichartz
} inequalities

\begin{lemma}[Strichartz estimate\cite{KeT98}, \cite{Stri77}]\label{se}
Let $u$ be an $\dot{S}^0$ solution to the Schr\"{o}dinger equation
$(\ref{equ1})$. Then
\begin{equation*}
\big\|u\big\|_{\dot{S}^0} \lesssim \big\|u(t_0)\big\|_{L^2_x} +
\big\|f(u)\big\|_{L^{q'}_tL^{r'}_x(I\times \mathbb{R}^d)}
\end{equation*}
for any $t_0 \in I$ and any admissible pairs $(q, r)$. The implicit
constant is independent of the choice of interval $I$.
\end{lemma}

From Sobolev embedding, we have

\begin{lemma}\label{et}
For any function $u$ on $I\times \mathbb{R}^d$, we have
\begin{equation*}
\big\|\nabla u\big\|_{L^{\infty}_tL^2_x} + \big\|\nabla
u\big\|_{L^6_tL^{\frac{6d}{3d-2}}_x}+ \big\|\nabla
u\big\|_{L^3_tL^{\frac{6d}{3d-4}}_x} + \big\|
u\big\|_{L^{\infty}_tL^{\frac{2d}{d-2}}_x} + \big\|
u\big\|_{L^6_tL^{\frac{6d}{3d-8}}_x} \lesssim \big\|\nabla
u\big\|_{\dot{S}^0},
\end{equation*}
where all spacetime norms are on $I\times \mathbb{R}^d$.
\end{lemma}

For convenience, we introduce two abbreviated notations. For a time
interval $I$, we denote
\begin{equation*}
\aligned
 \big\|u\big\|_{X(I)}  :=
\big\|u\big\|_{L^6_t(I;L^{\frac{6d}{3d-8}}_x )};\quad  \big\|u\big\|_{Y(I)}  :=
\big\|\nabla u\big\|_{L^6_t(I;L^{\frac{6d}{3d-2}}_x )}; \quad
 \big\|u\big\|_{W(I)} :=
\big\|\nabla u\big\|_{L^3_t(I;L^{\frac{6d}{3d-4}}_x )}.
\endaligned
\end{equation*}

We develop a local well-posedness and blow-up criterion for the
$\dot{H}^1$-critical Hartree equation. First, we have

\begin{proposition}[Local well-posedness \cite{MiXZ07a}]\label{lwp}
Suppose $\big\|u(t_0)\big\|_{\dot{H}^1} \leq A$, $I$ be a compact time interval that
contains $t_0$ such that
\begin{equation*}
\big\| U(t-t_0) u(t_0) \big\|_{X(I)}\leq \delta,
\end{equation*}
for a sufficiently small absolute constant $\delta=\delta(A)>0$. Then there exists a
unique solution $u \in C^0_t \dot{H}^1_x$ to $(\ref{equ1})$ on $I\times \mathbb{R}^d$,
such that
\begin{equation*}
\big\| u \big\|_{W(I)}< \infty, \quad \big\| u\big\|_{X(I)} \leq
2\delta.
\end{equation*}
Moreover, if $u_{0,k} \rightarrow u_0$ in $\dot{H}^1(\mathbb{R}^d)$, the
corresponding solutions $u_k \rightarrow u$ in $C\big(I;
\dot{H}^1(\mathbb{R}^d)\big)$.
\end{proposition}

\begin{remark}\label{sdw}
There exists $\widetilde{\delta}>0$, such that if
$\big\|u(t_0)\big\|_{\dot{H}^1} \leq \widetilde{\delta}$, the
conclusion of Proposition \ref{lwp} applies to any interval $I$. In
fact, by Strichartz estimates, we have
\begin{equation*}
\big\| e^{i(t-t_0)\Delta}u(t_0)\big\|_{X(I)} \leq C  \big\|
e^{i(t-t_0)\Delta}u(t_0)\big\|_{Y(I)} \leq C \widetilde{\delta},
\end{equation*}
and the claim follows.
\end{remark}

\begin{remark}\label{interexist}
Given $u_0 \in \dot{H}^1$, there exists $I$ such that $ 0\in I$ and the
hypothesis of Proposition \ref{lwp} is satisfied on $I$. In fact, by Strichartz estimates, we have
\begin{equation*}
\aligned
\big\|
  e^{it\Delta}u_0 \big\|_{Y(I)}
< \infty,
\endaligned
\end{equation*}
then the claim follows from Sobolev inequality and absolutely continuity theorem.
\end{remark}

\begin{remark}[Energy identity]\label{ei}
Based on the standard limiting argument, if $u$ is the solution
constructed in Proposition $\ref{lwp}$, we have that
\begin{equation}\label{energy}
\aligned E(u(t))=\frac12 & \big\|\nabla
u(t)\big\|^2_{L^2}-\frac{1}{4} \iint  \frac{1}{|x-y|^{4}}
|u(t,x)|^2 |u(t,y)|^2\ dxdy.
\endaligned
\end{equation}
is constant for $t\in I$.
\end{remark}

Now let $t_0\in I$. We say that $u\in C(I; \dot{H}^1(\mathbb{R}^d)) \cap
W(I)$ is a solution of $(\ref{equ1})$ if
\begin{equation*}
\aligned u\mid_{t_0}=u_0, \ \ \text{and}\ \
u(t)=e^{i(t-t_0)\Delta}u_0 -i \int^t_{t_0}e^{i(t-s)\Delta}f(u)ds
\endaligned
\end{equation*}
with $f(u)=-\big( |x|^{-4}* |u|^2 \big) u$. Note that if $u^{(1)}$,
$u^{(2)}$ are solutions of $(\ref{equ1})$ on $I$,
$u^{(1)}(t_0)=u^{(2)}(t_0)$, then $u^{(1)}\equiv u^{(2)}$ on
$I\times \mathbb{R}^d$. This is because we can partition $I$ into a
finite collection of subintervals $I_j$ with
\begin{equation*}
\aligned
 A=\sup_{t\in I}\max_{i=1,2}
\big\|u^{(i)}(t)\big\|_{\dot{H}^1}.\endaligned
\end{equation*}
%the $S(I_j)$ norm and the
%$W(I_j)$ norm  of $u^{(i)}$ are less than $a, b$, where $a, b$ are
%obtained in the proof of Proposition \ref{lwp}.
If $j_0$ is  such that $t_0 \in I_{j_0}$, then the uniqueness of the
fixed point in the proof of Proposition \ref{lwp}, combined with
Remark \ref{interexist} gives an interval $\widetilde{I} \ni t_0$ so
that $u^{(1)}(t) = u^{(2)}(t), t\in \widetilde{I}$. A continuation
argument now easily gives $u^{(1)}(t) = u^{(2)}(t), t\in I$.

\begin{definition}[Maximal interval]\label{mi}
The above analysis allows us to define a maximal interval $\big(t_0-T_{-}(u_0),t_0+
T_{+}(u_0)\big)$, with $T_{\pm}(u_0)>0$, where the solution is defined.
If $T_1<t_0+T_{+}(u_0)$, $T_2>t_0-T_{-}(u_0)$, $T_2<t_0<T_1$, then
$u$ solves $(\ref{equ1})$ in $[T_2, T_1]\times \mathbb{R}^d$, so
that $u\in C([T_2, T_1], \dot{H}^1(\mathbb{R}^d))\cap X([T_2, T_1])
\cap W([T_2, T_1])$.
\end{definition}

\begin{proposition}[Blow-up criterion \cite{MiXZ07a}]\label{buc}
If $T_{+}(u_0)<+\infty$, then
\begin{equation*}
\big\| u \big\|_{X\big(t_0, t_0+T_{+}(u_0)\big)} =+\infty.
\end{equation*}
A corresponding result holds for $T_{-}(u_0)$.
\end{proposition}

\begin{definition}[Nonlinear profile]\label{np}
Let $v_0 \in \dot{H}^1$, $v(t)=e^{it\Delta}v_0$ and let $t_n$ be a
sequence, with $\displaystyle \lim_{n\rightarrow \infty}t_n =
\overline{t} \in [-\infty, \infty]$. We say that $u(t,x)$ is a
nonlinear profile associated with $(v_0,\{t_n\})$ if there exists an
interval $I$, with $\overline{t} \in I$ (if $\overline{t}=\pm
\infty, I=[a, +\infty)$ or $(-\infty, a]$) such that $u$ is a
solution of $(\ref{equ1})$ in $I$ and
\begin{equation*}
\aligned \lim_{n\rightarrow \infty} \big\| u(t_n, \cdot)-v(t_n,
\cdot)\big\|_{\dot{H}^1} =0.
\endaligned
\end{equation*}
\end{definition}

\begin{remark}
Similar as in \cite{KeM06}, there always exists a unique nonlinear
profile $u(t)$ associated to $(v_0, \{t_n\})$, with a maximal
interval $I$.
\end{remark}

Last, in order to meet our needs in Lemma \ref{sc}, we give a stability theory,
which is somewhat different from that in \cite{MiXZ07b}, but their proofs are similar in essence.

\begin{proposition}[Long-time
perturbations ]\label{ltp} Let I be a compact
interval, and let $\widetilde{u}$ be a function on $I\times
\mathbb{R}^d$ which obeys the bounds
\begin{equation}\label{stc}
\big\| \widetilde{u} \big\|_{X(I)}   \leq M
\end{equation}
and
\begin{equation}\label{kc}
\big\| \widetilde{u} \big\|_{L^{\infty}_t( I;\dot{H}^1_x)} \leq E
\end{equation}
for some $M, E>0$. Suppose also that $\widetilde{u}$ is a
near-solution to $(\ref{equ1})$ in the sense that it solves
\begin{equation}\label{near-solution}
\aligned (i\partial_t + \Delta ) \widetilde{u} =-
(|x|^{-4}*|\widetilde{u}|^2)\widetilde{u}+e
\endaligned
\end{equation}
 for some function $e$. Let $t_0\in I$, and
let $u(t_0)$ be close to $\widetilde{u}(t_0)$ in the sense that
\begin{equation*}
\big\| u(t_0)-\widetilde{u}(t_0) \big\|_{\dot{H}^1_x}
\leq E'
\end{equation*}
for some $E'>0$. Assume also that we have the smallness conditions
\begin{eqnarray}
\big\| e^{i(t-t_0)\Delta}
 \big(u(t_0)-\widetilde{u}(t_0)\big)
\big\|_{{\cal Z}^1(I)} & \leq &
\epsilon , \label{longcond4}\\
 \big\| e \big\|_{L^{\frac32}_t(I;\dot{H}_x^{1, \frac{6d}{3d+4}}) }
& \leq  &\epsilon     \label{ec}
\end{eqnarray}
for some $0< \epsilon< \epsilon_1$, where $\epsilon_1$ is some
constant $\epsilon_1=\epsilon_1(E, E', M)>0$.

We conclude that there exists a solution $u$ to $(\ref{equ1})$ on
$I\times \mathbb{R}^d$ with the specified initial data $u(t_0)$ at
$t_0$, and
\begin{equation*}\aligned
\big\| u \big\|_{{\cal Z}^1(I)}  \leq C(M, E, E').
\endaligned\end{equation*}
Moreover, we have
\begin{equation*}
\aligned
\big\|\nabla u \big\|_{S^0(I)}  \leq C(M, E, E').
\endaligned
\end{equation*}
\end{proposition}

\begin{remark}\label{stab1}
Under the assumptions $(\ref{stc})$ and $(\ref{ec})$, we know that the assumption $(\ref{kc})$ is equivalent to the following condition
\begin{equation*}
\aligned
\big\|\nabla \widetilde{u}(t_0)\big\|_{L^2} \leq E.
\endaligned
\end{equation*}
\end{remark}

\begin{remark}\label{ct}
The long time perturbation theorem in \cite{MiXZ07b} yields the following continuity fact, which
will be used later: Let $\widetilde{u}_0\in \dot{H}^1$, $\big\|
\widetilde{u}_0 \big\|_{\dot{H}^1}\leq A,$ and let $\widetilde{u}$
be the solution of $(\ref{equ1})$, with maximal interval of
existence $\big( T_{-}(\widetilde{u}_0 ), T_{+}(\widetilde{u}_0)
\big)$. Let $u_{0,n}\rightarrow \widetilde{u}_0$ in $\dot{H}^1$, and
let $u_n$ be the corresponding solution of $(\ref{equ1})$, with maximal interval of existence $\big( T_{-}(u_{0,n} ), T_{+}(u_{0,n})
\big)$. Then
\begin{equation*}
\aligned T_{-}(\widetilde{u}_0 ) & \geq \overline{\lim_{n
\rightarrow +\infty}}
T_{-}(u_{0,n} ), \\
T_{+}(\widetilde{u}_0) & \leq \lim_{ \overline{n \rightarrow
+\infty}} T_{+}(u_{0,n} ),
\endaligned
\end{equation*}
and for each $t \in \big( T_{-}(\widetilde{u}_0 ),
T_{+}(\widetilde{u}_0) \big)$, $u_n(t) \rightarrow \widetilde{u}(t)$
in $\dot{H}^1$.
\end{remark}

%%%%%%%%%%%%%%%%%%%%%%%%%%%%%%%%%%%%%%%%%%%%%%%%%%%%%%%%%%%%%%%%%%%%%%%%%%%%%%%%%%%
\section{Some variational estimates and blow-up result}
\setcounter{section}{3} \setcounter{equation}{0} Let
$W(x)$ be the ground state to be the positive radial Schwartz solution to the elliptic equation
\begin{equation} \label{ellipe}
 \Delta W  + \big( |x|^{-4}* |W|^2 \big) W =0.
\end{equation}
The existence and uniqueness of $W$ were established in \cite{Lieb83} and \cite{Liu07}. By invariance of the equation, for
$\theta_0 \in [-\pi, \pi]$, $\lambda_0 > 0$, $x_0 \in \mathbb{R}^d$,
\begin{equation*}
W_{\theta_0, x_0,
\lambda_0}(x)=\lambda_0^{-\frac{d-2}{2}}e^{i\theta_0}W\big(\frac{x-x_0}{\lambda_0}\big)
\end{equation*}
is still a solution. Now let $C_d$ be the best constant of the Sobolev inequality in dimension
$d$. That is,
\begin{equation}\label{SobolevIne}
\forall u\in \dot{H}^1, \quad \big\| \big( |x|^{-4}* |u|^2 \big)
|u|^2\big\|^{\frac14}_{L^1} \leq C_d \big\| \nabla u\big\|_{L^2}.
\end{equation}
In addition, using the
concentration-compactness argument \cite{HmK05}, \cite{Lions84}, \cite{Lions85} and \cite{MiXZ07c}, we can obtain the
following characterization of $W$:

If $\big\| \big( |x|^{-4}* |u|^2 \big) |u|^2\big\|^{\frac14}_{L^1} =
C_d \big\| \nabla u\big\|_{L^2}, u\not =0$, then $\exists (\theta_0,
\lambda_0, x_0 )$ such that $u=W_{\theta_0, x_0, \lambda_0}$.

From above, we have
\begin{equation*}
\aligned \big\| \big( |x|^{-4}* |W|^2 \big) |W|^2\big\|_{L^1}&=C^4_d
\big(\int |\nabla W|^2 dx \big)^2.
\endaligned
\end{equation*}
On the other hand,  from $(\ref{ellipe})$, we obtain
\begin{equation*}
\aligned \big\| \big( |x|^{-4}* |W|^2 \big) |W|^2\big\|_{L^1}&=\int
|\nabla W|^2 dx.
\endaligned
\end{equation*}
Hence, we have
\begin{equation*}
\aligned \big\| \nabla W \big\|^2_{L^2} = \frac{1}{C^{4}_d},\quad
E(W)= (\frac12 - \frac14) \big\| \nabla W \big\|^2_{L^2} =
\frac{1}{4 C^4_d}.
\endaligned
\end{equation*}

\begin{lemma}\label{slbv}
Assume that
\begin{equation*}
\big\|\nabla u \big\|_{L^2} < \big\|\nabla W \big\|_{L^2}.
\end{equation*}
Assume moreover that $E(u) \leq (1-\delta_0)E(W)$ where
$\delta_0>0$. Then, there exists
$\overline{\delta}=\delta^{1/2}_0>0$ such that
\begin{equation*}
\aligned \int |\nabla u |^2dx  - \iint &
\frac{|u(x)|^2|u(y)|^2}{|x-y|^4} dxdy \geq \frac{\overline{\delta}}{2} \int
|\nabla u |^2dx, \\
\int |\nabla u |^2dx & \leq (1-
\overline{\delta})
\int  |\nabla W  |^2dx ,\\
E(u) &  \geq 0.
\endaligned
\end{equation*}
\end{lemma}

{\bf Proof: }
Define
\begin{equation*}
\aligned a= \int |\nabla u |^2dx \quad \text{and}\ \
f(x)=\frac12x-\frac14C^4_d x^2.
\endaligned
\end{equation*}
From (\ref{SobolevIne}), we have
\begin{equation}\label{constrain}
\aligned
(1-\delta_0)E(W) \geq  E(u) \geq \frac12\int |\nabla u |^2dx
-\frac14C^4_d \big( \int |\nabla u |^2dx \big)^2=f(a).
\endaligned
\end{equation}

Note that
\begin{equation*}
\aligned f'(x)=\frac12-\frac12C^4_d x,
\endaligned
\end{equation*}
This implies that
\begin{equation*}
\aligned f'(x)=0 \Longleftrightarrow x=\frac{1}{C^4_d}=\int |\nabla
W(x)|^2dx.
\endaligned
\end{equation*}
On the other hand,
\begin{equation*}
\aligned f'(x)> 0,  \ &\text{for}\ \ x< \frac{1}{C^4_d}, \\
f(0)=0, \ \ f\big( \frac{1}{C^4_d} \big)& = \frac{1}{4C^4_d}=E(W).
\endaligned
\end{equation*}
Together with (\ref{constrain}) and the fact that $a=\big\|\nabla u\big\|^2_{L^2} \in \big[0, \frac{1}{C^4_d}\big)$, these imply that
\begin{equation*}
\aligned \big\|\nabla u\big\|^2_{L^2}=a  \leq (1-&
\overline{\delta}) \frac{1}{C^4_d}=(1-
\overline{\delta})
\int  |\nabla W  |^2dx, \ \overline{\delta} = \delta^{1/2}_0,\\
E(u)& \geq f(a)\geq 0.
\endaligned
\end{equation*}

Now define
\begin{equation*}
\aligned g(x)=x-C^4_dx^2.
\endaligned
\end{equation*}
From (\ref{SobolevIne}), we also have
\begin{equation}\label{constrain2}
\aligned
\int |\nabla u |^2dx  - \iint &
\frac{|u(x)|^2|u(y)|^2}{|x-y|^4} dxdy \geq \int |\nabla u |^2dx  - C^4_d \Big( \int |\nabla u |^2dx  \Big)^2 =g\big(a\big).
\endaligned
\end{equation}

Note that
\begin{equation*}
\aligned
g(x)=0 \Longleftrightarrow x=0,\ & \text{or}\ x=\frac{1}{C^4_d},  \\
g'(0)=1,\   g'(\frac{1}{C^4_d})=-1, \ & g''(x)=-2C^4_d<0.
\endaligned
\end{equation*}
Hence, we obtain
\begin{equation*}
\aligned
g(x)\geq \frac12 \min\big(x, \frac{1}{C^4_d} -x \big) \quad \text{for} \ \ 0\leq x\leq \frac{1}{C^4_d}.
\endaligned
\end{equation*}
Since $\big\|\nabla u\big\|^2_{L^2} = a \in [0,(1-
\overline{\delta}) \frac{1}{C^4_d}]$, the above inequality implies that
\begin{equation*}
\aligned
\text{(LHS)}\ of (\ref{constrain2}) \geq g(a) & \geq \frac12 \min(a, \frac{1}{C^4_d} -a) \\
&\geq \frac12 \min(a, \overline{\delta} a) = \frac{\overline{\delta}}{2}  a.
\endaligned
\end{equation*}
This completes the proof.

\begin{corollary}\label{pe}
Assume that $u\in \dot{H}^1(\mathbb{R}^d)$ and that $\big\|\nabla u \big\|_{L^2} <
\big\|\nabla W \big\|_{L^2}$. Then $E(u)\geq 0$.
\end{corollary}

{\bf Proof: } If $E(u)<E(W)$, the conclusion follows from Lemma \ref{slbv}.
If $E(u)\geq E(W)=\frac{1}{4C^4_d}$, it is clear.

\begin{proposition}[Lower bound on the convexity of the variance]\label{dlbv}
Let $u$ be a solution of $(\ref{equ1})$ with $t_0=0, u(0)=u_0$ such
that for $\delta_0 >0$
\begin{equation*}
\int |\nabla u_0 |^2dx < \int  |\nabla W |^2dx, \quad
E(u_0)<(1-\delta_0)E(W).
\end{equation*}
Let $I \ni 0$ be the maximal interval of existence given by
Definition \ref{mi}. Let $\overline{\delta}=\delta^{1/2}_0$
be as in Lemma \ref{slbv}. Then for each $t\in I$, we have
\begin{equation*}
\aligned \int |\nabla u (t) |^2dx  - \iint &
\frac{|u(t)|^2|u(t)|^2}{|x-y|^4} dxdy \geq \frac{\overline{\delta}}{2} \int
|\nabla u(t) |^2dx, \\
 \int |\nabla u (t)|^2dx & \leq (1-
\overline{\delta})
\int  |\nabla W  |^2dx ,\\
E(u(t)) &  \geq 0.
\endaligned
\end{equation*}
\end{proposition}

{\bf Proof: } We prove it by the continuity argument.
Define
\begin{equation*}
\aligned
\Omega = \big\{t \in I, \big\|\nabla u(t)\big\|_{L^2} < \big\| \nabla W \big\|_{L^2}, E(u(t))< (1-\delta_0)E(W) \big\}.
\endaligned
\end{equation*}
It suffices to prove that $\Omega$ is both open and closed.

Firstly, we see that $t_0 \in \Omega$. Secondly, $\Omega$ is open because of $u\in C^0_t(I, \dot{H}^1)$ and the conservation of energy.
Lastly, we need to prove that $\Omega$ is also closed. For any $t_n \in \Omega, T\in I$, and $t_n \rightarrow T$. Then
\begin{equation*}
\aligned
\big\| \nabla u(t_n) \big\|_{ L^2} < \big\| \nabla W\big\|_{L^2},\quad
E(u(t_n))  < (1-\delta_0)E(W).
\endaligned
\end{equation*}
From Lemma \ref{slbv}, we obtain
\begin{equation*}
\aligned
\big\| \nabla u(t_n) \big\|^2_{L^2} & < (1-\overline{\delta}) \big\| \nabla W\big\|^2_{L^2}.
\endaligned
\end{equation*}
Using the fact that $u\in C^0_t(I, \dot{H}^1)$ and the conservation of energy again, we have
\begin{equation*}
\aligned
\big\| \nabla u(T) \big\|^2_{L^2} & \leq (1-\overline{\delta}) \big\| \nabla W\big\|^2_{L^2},\quad E(u(T))=E(u(t_n))  < (1-\delta_0)E(W).
\endaligned
\end{equation*}
This implies that $T\in \Omega$ and completes the proof.

\begin{corollary}[Comparability of gradient and energy]\label{cge}
Let $u, u_0$ be as in Proposition \ref{dlbv}. Then for all $t\in I$
we have
\begin{equation*}
E(u(t)) \thickapprox \int | \nabla u(t)|^2dx \thickapprox \int |
\nabla u_0|^2dx
\end{equation*}
with comparability constants which depend only on $\delta_0$.
\end{corollary}

{\bf Proof: } From Proposition \ref{dlbv}, we have
\begin{equation*}
\aligned
\frac12\int |\nabla u (t) |^2dx \geq E(u(t))& =\frac14 \int |\nabla u (t) |^2dx +\frac14\big( \int |\nabla u (t,x) |^2dx  - \iint &
\frac{|u(t,x)|^2|u(t,y)|^2}{|x-y|^4} dxdy \big) \\
&  \geq \frac{2+\overline{\delta}}{8} \int
|\nabla u(t) |^2dx \quad \forall\ t\in I.
\endaligned
\end{equation*}
This together with the conservation of energy implies the claim.

In order to obtain blow up results, we first give the (local) virial identity, which we can verify by some direct computations.
\begin{lemma}\label{virialidentity}
Let $ \varphi\in C^{\infty}_{0} (\mathbb{R}^d)$, $V(x)=|x|^{-4}$, $t \in [0,
T_{+}(u_0))$. Then
\begin{equation*}
\aligned
(1)& \quad \frac{d}{dt}\int \big| u \big|^2 \varphi dx \ = 2 \text{Im}
\int \overline{u} \nabla u \nabla\varphi dx;\\
(2)& \quad \frac{d^2}{dt^2}\int \big| u \big|^2 \varphi dx
=-\int\triangle\triangle \varphi |u|^{2}dx+4\mathrm{Re}\int
\varphi_{jk}\overline{u}_{j}u_{k}dx \\
& \qquad \qquad \qquad \qquad \ -\mathrm{Re}\int\int \big(\nabla \varphi(x)- \nabla
\varphi(y)\big)\nabla
V(x-y)|u(y)|^2|u(x)|^{2}dxdy.
\endaligned
\end{equation*}
\end{lemma}

\begin{proposition}\label{blowupresults}
Assume that $u_0 \in \dot{H}^1(\mathbb{R}^d)$ and
\begin{equation*}
\aligned
 E(u_0)<E(W),\quad \int |\nabla u_0 |^2dx > \int  |\nabla W
|^2dx.
\endaligned
\end{equation*}
If $|x|u_0 \in L^2$ or $u_0$ is radial, then the maximal
interval I of existence  must be finite.
\end{proposition}

{\bf Proof: } Indeed, we can choose a suitable small number $\delta_0>0$, such that
\begin{equation*}
\aligned
E(u_0)<(1-\delta_0)E(W),\quad \int |\nabla u_0 |^2dx > \int  |\nabla W
|^2dx.
\endaligned
\end{equation*}
Arguing as in Lemma \ref{slbv}, we obtain that there exists $\widetilde{\delta}$ such that
\begin{equation*}
\aligned
\int |\nabla u_0 |^2dx > (1+\widetilde{\delta})\int  |\nabla W
|^2dx = \frac{1+\widetilde{\delta}}{C^4_d}.
\endaligned
\end{equation*}
This shows that
\begin{equation*}
\aligned
\int |\nabla u_0|^2dx  - \iint
\frac{|u_0(x)|^2|u_0(y)|^2}{|x-y|^4} dxdy & = 4E(u_0)-\int |\nabla u_0|^2dx \qquad \qquad \qquad  \\
& < 4(1-\delta_0)E(W) - \frac{1+\widetilde{\delta}}{C^4_d}  = \frac{1-\delta_0}{C^4_d} - \frac{1+\widetilde{\delta}}{C^4_d} \\
& = - \frac{\delta_0+\widetilde{\delta}}{C^4_d}<0.
\endaligned
\end{equation*}

Now define
\begin{equation*}
\aligned
\Omega = \big\{t \in I, \big\|\nabla u(t)\big\|_{L^2} > \big\| \nabla W \big\|_{L^2}, E(u(t))< (1-\delta_0)E(W) \big\}.
\endaligned
\end{equation*}
Using the continuity argument and arguing as in Proposition \ref{dlbv}, we have
\begin{equation*}
\aligned
\Omega =I.
\endaligned
\end{equation*}
Arguing as in Lemma \ref{slbv} again, we have
\begin{equation*}
\aligned
\big\| \nabla u(t) \big\|^2_{L^2} > (1+\widetilde{\delta}) \big\| \nabla W\big\|^2_{L^2}.
\endaligned
\end{equation*}
Then
\begin{equation*}
\aligned
\int |\nabla u(t,x)|^2dx  - \iint
\frac{|u(t,x)|^2|u(t,y)|^2}{|x-y|^4} dxdy  = - \frac{\delta_0+\widetilde{\delta}}{C^4_d}<0, \forall\ t\in I.
\endaligned
\end{equation*}

As for the case that $|x|u_0 \in L^2$. From Lemma \ref{virialidentity}, we have
\begin{equation*}
\aligned
\frac{d^2}{dt^2} \int |x|^2 |u(t,x)|^2 dx = 8 \Big( \int |\nabla u(t,x)|^2dx  - \iint
\frac{|u(t,x)|^2|u(t,y)|^2}{|x-y|^4} dxdy  \Big) < 0.
\endaligned
\end{equation*}
This implies that $I$ must be finite.

As for the case that  $u_0 $ is radial. Using the local virial identity \cite{DuHR07}, \cite{HoR07} and \cite{OT91}, we can also deduce the same result.

%%%%%%%%%%%%%%%%%%%%%%%%%%%%%%%%%%%%%%%%%%%%%%%%%%%%%%%%%%%%%%%%%%%%%%%%%%%%%%%%%%%
\section{Existence and compactness of a critical element}
\setcounter{section}{4} \setcounter{equation}{0}

Let us consider the statement

\begin{enumerate}
\item[(SC)] For all $u_0 \in \dot{H}^1(\mathbb{R}^d)$ with $\big\|\nabla u_0 \big\|_{L^2} < \big\|\nabla W \big\|_{L^2},
E(u_0)< E(W)$, if $u$ is the corresponding solution to (\ref{equ1}),
with maximal interval of existence $I$, then $I=(-\infty, +\infty)$
and $\big\| u \big\|_{X(\mathbb{R})}<+\infty$.
\end{enumerate}

We say that $(SC)(u_0)$ holds if for this particular $u_0$ with
$\big\|\nabla u_0 \big\|_{L^2} < \big\|\nabla W \big\|_{L^2}$, $E(u_0)< E(W)$, and
$u$ is the corresponding solution to (\ref{equ1}), with maximal
interval of existence $I$, then $I=(-\infty, +\infty)$ and $\big\| u
\big\|_{X(\mathbb{R})}<+\infty$.

Note that, because of Remark \ref{sdw}, if $\big\| u_0
\big\|_{\dot{H}^1}\leq \widetilde{\delta}$,  $(SC)(u_0)$ holds.
Thus, in light of Corollary \ref{cge}, there exists $\eta_0 > 0$
such that if $u_0$ is as in $(SC)$ and $E(u_0) < \eta_0$, then
$(SC)(u_0)$ holds. Moreover, $E(u_0)\geq
0$ in light of Proposition \ref{dlbv}. Thus, there exists a number
$E_c$, with $\eta_0 \leq E_c \leq E(W)$, such that, if $u_0$ is
radial with $\big\|\nabla u_0 \big\|_{L^2} < \big\|\nabla W \big\|_{L^2},
E(u_0)< E_c$, then $(SC)(u_0)$ holds, and $E_c$ is optimal with this
property. If $ E_c \geq  E(W)$, then the first part of Theorem \ref{main} is true. For the
rest of this section, we will assume that $E_c < E(W)$ and
ultimately deduce a contradiction in Section 5. By definition of
$E_c$, we have
\begin{enumerate}
\item[(C.1)]If $u_0$ is radial and $\big\|\nabla u_0 \big\|_{L^2} < \big\|\nabla W \big\|_{L^2},
E(u_0)< E_c$, then $(SC)(u_0)$ holds.
\item[(C.2)]There exists a sequence of radial solutions $u_n$ to
(\ref{equ1}) with corresponding initial data $u_{n,0}$ such that
$\big\|\nabla u_{n,0} \big\|_{L^2} < \big\|\nabla W \big\|_{L^2},
E(u_{n,0})\searrow E_c$ as $n \rightarrow +\infty$, for which
$(SC)(u_{n,0})$ does not  hold for any $n$.
\end{enumerate}

The goal of this section is to use the above sequence $u_{n,0}$ to
prove the existence of an $\dot{H}^1$ radial solution $u_c$ to
(\ref{equ1}) with initial data $u_{c,0}$ such that $\big\|\nabla
u_{c,0} \big\|_{L^2} < \big\|\nabla W \big\|_{L^2}$, $E(u_{c,0}) = E_c$ for which
$(SC)(u_{c,0})$ does not  hold (see Proposition \ref{ecs}). Moreover, we will show that this
critical solution has a compactness property up to the symmetries of
this equation (see Proposition \ref{ccs}).

Before stating and proving Proposition \ref{ecs}, we introduce some useful
preliminaries in the spirit of the results of Keraani \cite{Ker01}. First we give the profile decomposition lemma.

\begin{lemma}[Profile decomposition]\label{pd}
Let $v_{n,0}$ be a radial uniformly bounded sequence in $\dot{H}^1$,
i.e. $\big\|\nabla v_{n,0}\big\|_{L^2} \leq A$. Assume that
$\big\|e^{it\Delta}v_{n,0}\big\|_{X(\mathbb{R})}\geq \delta>0$,
where $\delta=\delta(d)$ is as in Proposition \ref{lwp}. Then for
each $J$, there exists a subsequence of $v_{n,0}$, also denoted
$v_{n,0}$, and
\begin{enumerate}
\item[$(1)$]For each $1\leq j\leq J$, there exists a radial profile $  V_{0,j } $ in $\dot{H}^1$.
\item[$(2)$]For each $1\leq j\leq J$, there exists a sequence of $(\lambda_{j,n},
t_{j,n})$ with
\begin{equation}\label{orthogonality}
\aligned
\frac{\lambda_{j,n}}{\lambda_{j',n}} +
\frac{\lambda_{j',n}}{\lambda_{j,n}} + \frac{|t_{j,n}
-t_{j',n}|}{\lambda^2_{j,n}} \rightarrow \infty \quad \text{as} \ \
n \rightarrow \infty \quad \text{for} \ \ j\not=j'.
\endaligned
\end{equation}
\item[$(3)$]There exists a sequence of radial remainder $w^J_n$ in
$\dot{H}^1$,
\end{enumerate}
such that
\begin{equation}\label{expansion}
\aligned
v_{n,0}(x)&=\sum^{J}_{j=1}\frac{1}{\lambda^{(d-2)/2}_{j,n}}
e^{-it_{j,n}\Delta}V_{0,j}(\frac{x}{\lambda_{j,n}} \big) + w^J_n(x)\\
&=\sum^{J}_{j=1}\frac{1}{\lambda^{(d-2)/2}_{j,n}}
V^{l}_{j}\big(-\frac{t_{j,n}}{\lambda^2_{j,n}},
\frac{x}{\lambda_{j,n}} \big) + w^J_n(x)
\endaligned
\end{equation}
with
\begin{eqnarray}
V^l_j(t,x) =e^{it\Delta }V_{0,j}(x), & &  \big\|V_{0,1}\big\|_{\dot{H}^1}  \geq \alpha_0(A)>0, \label{fp}\\
\big\|\nabla v_{n,0} \big\|^2_{L^2} & = &\sum^J_{j=1} \big\|\nabla
V_{0,j}
\big\|^2_{L^2} + \big\|\nabla w^J_n \big\|^2_{L^2} +o_n(1), \label{pathak}\\
E(v_{n,0}) & = & \sum^J_{j=1}
E(V^l_j(-\frac{t_{j,n}}{\lambda^2_{j,n}})) + E(w^J_n) +o_n(1), \label{pathae}\\
\lim_{J \rightarrow \infty}  \big[ \lim_{n \rightarrow \infty}  \big\|
e^{it\Delta} w^J_n \big\|_{L^q(\mathbb{R},L^r)} \big]  &=&0, \ \forall \  \frac2q = d(\frac12-\frac1r)-1, \ \frac{2d}{d-2}\leq r < \frac{2d}{d-4}. \label{asym}
\end{eqnarray}
\end{lemma}

{\bf Proof: } Here we only give the proof of energy asymptotic Pythagorean expansion (\ref{pathae}), the rest is standard (see \cite{Ker01}).

By the asymptotic Pythagorean expansion of kinetic energy, it suffices to show that
\begin{equation*}
\aligned
\iint  \frac{1}{|x-y|^{4}}
|v_{n,0}(x)|^2 |v_{n,0}(y)|^2\ dxdy & = \sum^J_{j=1} \iint \frac{1}{|x-y|^{4}}
|V^l_j(-\frac{t_{j,n}}{\lambda^2_{j,n}},x)|^2 |V^l_j(-\frac{t_{j,n}}{\lambda^2_{j,n}},y)|^2\ dxdy \\
&  + \iint  \frac{1}{|x-y|^{4}}
|w^J_n(x)|^2 |w^J_n(y)|^2\ dxdy +o_n(1), \ \forall \ J\geq 1.
\endaligned
\end{equation*}

We first claim that if $J\geq 1$ is fixed, the orthogonality condition (\ref{orthogonality}) implies that
\begin{equation}\label{claim1}
\aligned
& \iint \frac{1}{|x-y|^{4}}
\Big|\sum^{J}_{j=1}\frac{1}{\lambda^{\frac{d-2}{2}}_{j,n}}
V^{l}_{j}\big(-\frac{t_{j,n}}{\lambda^2_{j,n}},
\frac{x}{\lambda_{j,n}} \big)\Big|^2 \Big|\sum^{J}_{j=1}\frac{1}{\lambda^{\frac{d-2}{2}}_{j,n}}
V^{l}_{j}\big(-\frac{t_{j,n}}{\lambda^2_{j,n}},
\frac{y}{\lambda_{j,n}} \big)\Big|^2\ dxdy \\
 = &\sum^J_{j=1} \iint \frac{1}{|x-y|^{4}}
|V^l_j(-\frac{t_{j,n}}{\lambda^2_{j,n}},x)|^2 |V^l_j(-\frac{t_{j,n}}{\lambda^2_{j,n}},y)|^2\ dxdy
 +o_n(1).
\endaligned
\end{equation}

By reindexing, we can arrange such that there is $J_0 \leq J$ with
\begin{enumerate}
\item[$(1)$] $ \forall\ 1 \leq j \leq J_0$, we have that $\displaystyle \Big| \frac{t_{j,n}}{\lambda^2_{j,n}} \Big| \leq C$ in $n$;
\item[$(2)$] $ \forall\ J_0 + 1 \leq j \leq J$, we have that $\displaystyle \Big| \frac{t_{j,n}}{\lambda^2_{j,n}} \Big|  \rightarrow +\infty $ as $n \rightarrow +\infty$.
\end{enumerate}
By passing to a subsequence and adjusting the profile $V_{0,j}$, we may assume that
\begin{equation*}
\aligned
 \forall \  1 \leq j \leq J_0, \quad \frac{t_{j,n}}{\lambda^2_{j,n}} =0,
\endaligned
\end{equation*}

From case (2), we have
\begin{equation}\label{high}
\aligned
\lim_{n\rightarrow +\infty}\iint \frac{1}{|x-y|^{4}}
|V^l_j(-\frac{t_{j,n}}{\lambda^2_{j,n}},x)|^2 |V^l_j(-\frac{t_{j,n}}{\lambda^2_{j,n}},y)|^2\ dxdy =0,\  \forall\ J_0 + 1 \leq j \leq J.
\endaligned
\end{equation}
Indeed, using Hardy inequality and the decay estimates for the free Schr\"{o}dinger equation
(similar to Lemma $4.1$ in \cite{GiV00} and Corollary 2.3.7 in \cite{Ca03}), we have for $J_0 + 1 \leq j \leq J$
\begin{equation*}
\aligned
\iint \frac{1}{|x-y|^{4}}
|V^l_j(-\frac{t_{j,n}}{\lambda^2_{j,n}},x)|^2 |V^l_j(-\frac{t_{j,n}}{\lambda^2_{j,n}},y)|^2\ dxdy &
\lesssim \big\| V^l_j(-\frac{t_{j,n}}{\lambda^2_{j,n}}) \big\|^4_{L^{\frac{2d}{d-2}}} \rightarrow 0, \ n\rightarrow +\infty.
\endaligned
\end{equation*}

By (\ref{orthogonality}), if $1\leq j < k \leq J_0$, we have
\begin{equation}\label{orthogonality2}
\aligned
\frac{\lambda_{j,n}}{\lambda_{k,n}} +
\frac{\lambda_{k,n}}{\lambda_{j,n}}  \rightarrow \infty \quad \text{as} \ \
n \rightarrow \infty.
\endaligned
\end{equation}
This implies that
\begin{equation}\label{low}
\aligned
& \iint \frac{1}{|x-y|^{4}}
\Big|\sum^{J_0}_{j=1}\frac{1}{\lambda^{(d-2)/2}_{j,n}}
V_{0,j}\big(\frac{x}{\lambda_{j,n}} \big)\Big|^2 \Big|\sum^{J_0}_{j=1}\frac{1}{\lambda^{(d-2)/2}_{j,n}}
V_{0,j}\big(\frac{y}{\lambda_{j,n}} \big)\Big|^2\ dxdy \\
 = &\sum^{J_0}_{j=1} \iint \frac{1}{|x-y|^{4}}
|V_{0,j}(x)|^2 |V_{0,j}(y)|^2\ dxdy
 +o_n(1).
\endaligned
\end{equation}
Hence, from (\ref{high}) and (\ref{low}), we obtain
\begin{equation*}
\aligned
& \iint \frac{1}{|x-y|^{4}}
\Big|\sum^{J}_{j=1}\frac{1}{\lambda^{\frac{d-2}{2}}_{j,n}}
V^{l}_{j}\big(-\frac{t_{j,n}}{\lambda^2_{j,n}},
\frac{x}{\lambda_{j,n}} \big)\Big|^2 \Big|\sum^{J}_{j=1}\frac{1}{\lambda^{\frac{d-2}{2}}_{j,n}}
V^{l}_{j}\big(-\frac{t_{j,n}}{\lambda^2_{j,n}},
\frac{y}{\lambda_{j,n}} \big)\Big|^2\ dxdy \\
 = & \iint \frac{1}{|x-y|^{4}}
\Big|\sum^{J_0}_{j=1}\frac{1}{\lambda^{\frac{d-2}{2}}_{j,n}}
V^{l}_{j}\big(-\frac{t_{j,n}}{\lambda^2_{j,n}},
\frac{x}{\lambda_{j,n}} \big)+ \sum^{J}_{j=J_0+1}\frac{1}{\lambda^{\frac{d-2}{2}}_{j,n}}
V^{l}_{j}\big(-\frac{t_{j,n}}{\lambda^2_{j,n}},
\frac{x}{\lambda_{j,n}} \big)\Big|^2 \\
& \qquad \qquad \ \times \Big|\sum^{J_0}_{j=1}\frac{1}{\lambda^{\frac{d-2}{2}}_{j,n}}
V^{l}_{j}\big(-\frac{t_{j,n}}{\lambda^2_{j,n}},
\frac{y}{\lambda_{j,n}} \big) + \sum^{J}_{j=J_0+1}\frac{1}{\lambda^{\frac{d-2}{2}}_{j,n}}
V^{l}_{j}\big(-\frac{t_{j,n}}{\lambda^2_{j,n}},
\frac{y}{\lambda_{j,n}} \big)\Big|^2\ dxdy \\
 = & \iint \frac{1}{|x-y|^{4}}
\Big|\sum^{J_0}_{j=1}\frac{1}{\lambda^{\frac{d-2}{2}}_{j,n}}
V_{0,j}\big(\frac{x}{\lambda_{j,n}} \big)+ \sum^{J}_{j=J_0+1}\frac{1}{\lambda^{\frac{d-2}{2}}_{j,n}}
V^{l}_{j}\big(-\frac{t_{j,n}}{\lambda^2_{j,n}},
\frac{x}{\lambda_{j,n}} \big)\Big|^2 \\
& \qquad \qquad \ \times \Big|\sum^{J_0}_{j=1}\frac{1}{\lambda^{\frac{d-2}{2}}_{j,n}}
V_{0,j}\big(\frac{y}{\lambda_{j,n}} \big) + \sum^{J}_{j=J_0+1}\frac{1}{\lambda^{\frac{d-2}{2}}_{j,n}}
V^{l}_{j}\big(-\frac{t_{j,n}}{\lambda^2_{j,n}},
\frac{y}{\lambda_{j,n}} \big)\Big|^2\ dxdy \\
=& \iint \frac{1}{|x-y|^{4}}
\Big|\sum^{J_0}_{j=1}\frac{1}{\lambda^{(d-2)/2}_{j,n}}
V_{0,j}\big(\frac{x}{\lambda_{j,n}} \big)\Big|^2 \Big|\sum^{J_0}_{j=1}\frac{1}{\lambda^{(d-2)/2}_{j,n}}
V_{0,j}\big(\frac{y}{\lambda_{j,n}} \big)\Big|^2\ dxdy \\
& \qquad \qquad \ +\sum^J_{j=J_0+1} \iint \frac{1}{|x-y|^{4}}
|V^l_j(-\frac{t_{j,n}}{\lambda^2_{j,n}},x)|^2 |V^l_j(-\frac{t_{j,n}}{\lambda^2_{j,n}},y)|^2\ dxdy +o_n(1)\\
=& \sum^J_{j=1} \iint \frac{1}{|x-y|^{4}}
|V^l_j(-\frac{t_{j,n}}{\lambda^2_{j,n}},x)|^2 |V^l_j(-\frac{t_{j,n}}{\lambda^2_{j,n}},y)|^2\ dxdy
 +o_n(1),
\endaligned
\end{equation*}
this yields (\ref{claim1}).

Secondly, we claim that
\begin{equation}\label{claim2}
\aligned
%\lim_{n\rightarrow +\infty} \iint  \frac{1}{|x-y|^{4}}
%|w^J_n(x)|^2 |w^J_n(y)|^2\ dxdy =
\lim_{n\rightarrow +\infty} \big\| w^J_n(x)\big\|_{L^{\frac{2d}{d-2}}_x} =0 \ \ \text{as} \ J\rightarrow +\infty.
\endaligned
\end{equation}
Indeed, we have
\begin{equation*}
\aligned
%\iint  \frac{1}{|x-y|^{4}}
%|w^J_n(x)|^2 |w^J_n(y)|^2\ dxdy & \lesssim
 \big\| w^J_n(x)\big\|_{L^{\frac{2d}{d-2}}_x} \lesssim  \big\| e^{it\Delta } w^J_n(x)\big\|_{L^{\infty}_t\big( \mathbb{R}; L^{\frac{2d}{d-2}}_x\big)},
\endaligned
\end{equation*}
this together with (\ref{asym}) implies the claim.

Note that (\ref{claim2}) implies that $\{w^J_n\}$ is uniformly bounded in $L^{\frac{2d}{d-2}}(\mathbb{R}^d)$,
the uniform boundness of $\{v_{n,0}\}$ in $\dot{H}^1(\mathbb{R}^d)$  also implies uniformly bounded  in $L^{\frac{2d}{d-2}}(\mathbb{R}^d)$.
Thus we can choose $J_1 \geq J$ and $N_1$ such that for $n\geq N_1$, we have
\begin{equation}\label{om1}
\aligned
& \ \Big| \iint  \frac{|v_{n,0}(x)|^2 |v_{n,0}(y)|^2}{|x-y|^{4}}
\ dxdy - \iint  \frac{|v_{n,0}(x)-w^{J_1}_n(x)|^2 |v_{n,0}(y)-w^{J_1}_n(y)|^2}{|x-y|^{4}}
\ dxdy \Big| \\
& + \Big| \iint  \frac{|w^J_n(x)-w^{J_1}_n(x)|^2 |w^J_n(y)-w^{J_1}_n(y)|^2}{|x-y|^{4}} \ dxdy -\iint  \frac{ |w^J_n(x)|^2 |w^J_n(y)|^2}{|x-y|^{4}}\ dxdy  \Big|\\
& \leq  C\big( \sup_{n} \big\| v_{n,0}(x) \big\|^3_{L^{\frac{2d}{d-2}}} + \sup_{n} \big\|w^{J}_n(x) \big\|^3_{L^{\frac{2d}{d-2}}}  \big)\big\|w^{J_1}_n(x)\big\|_{L^{\frac{2d}{d-2}}}
+ C \big\|w^{J_1}_n(x) \big\|^4_{L^{\frac{2d}{d-2}}} \leq  \epsilon.
\endaligned
\end{equation}
\ \ \ \ By (\ref{claim1}), we get $N_2 \geq N_1$ such that for $n\geq N_2$
\begin{equation}\label{om2}
\aligned
\Big| \iint & \frac{|v_{n,0}(x)-w^{J_1}_n(x)|^2 |v_{n,0}(y)-w^{J_1}_n(y)|^2}{|x-y|^{4}}
  dxdy \\
  & \qquad \qquad \qquad -\sum^{J_1}_{j=1} \iint \frac{|V^l_j(-\frac{t_{j,n}}{\lambda^2_{j,n}},x)|^2 |V^l_j(-\frac{t_{j,n}}{\lambda^2_{j,n}},y)|^2}{|x-y|^{4}}
 dxdy \Big| \leq \epsilon.
\endaligned
\end{equation}
\ \ \ \ Using (\ref{expansion}), we have
\begin{equation*}
\aligned
w^{J}_n(x)-w^{J_1}_n(x)= \sum^{J_1}_{j=J+1}\frac{1}{\lambda^{\frac{d-2}{2}}_{j,n}}
V^{l}_{j}\big(-\frac{t_{j,n}}{\lambda^2_{j,n}},
\frac{x}{\lambda_{j,n}} \big).
\endaligned
\end{equation*}
By (\ref{claim1}), we get $N_3 \geq N_2$ such that for $n\geq N_3$
\begin{equation*}\label{om3}
\aligned
\Big|\iint  \frac{|w^J_n(x)-w^{J_1}_n(x)|^2 |w^J_n(y)-w^{J_1}_n(y)|^2}{|x-y|^{4}} \ dxdy  - \sum^{J_1}_{j=J+1} \iint \frac{|V^l_j(-\frac{t_{j,n}}{\lambda^2_{j,n}},x)|^2 |V^l_j(-\frac{t_{j,n}}{\lambda^2_{j,n}},y)|^2}{|x-y|^{4}}
 dxdy  \Big| \leq \epsilon.
\endaligned
\end{equation*}
Combining the above inequality with (\ref{om1}), (\ref{om2}), we obtain that for $n\geq N_3$
\begin{equation*}
\aligned
 & \Big| \iint   \frac{|v_{n,0}(x)|^2 |v_{n,0}(y)|^2}{|x-y|^{4}}
 dxdy - \sum^J_{j=1} \iint \frac{|V^l_j(-\frac{t_{j,n}}{\lambda^2_{j,n}},x)|^2 |V^l_j(-\frac{t_{j,n}}{\lambda^2_{j,n}},y)|^2}{|x-y|^{4}}
 dxdy \\
 & \qquad \qquad \qquad  \qquad \qquad \qquad \ \ - \iint  \frac{|w^J_n(x)|^2 |w^J_n(y)|^2}{|x-y|^{4}} \ dxdy \Big| \\
=& \Big| \iint  \frac{|v_{n,0}(x)|^2 |v_{n,0}(y)|^2}{|x-y|^{4}}
 dxdy - \iint  \frac{|v_{n,0}(x)-w^{J_1}_n(x)|^2 |v_{n,0}(y)-w^{J_1}_n(y)|^2}{|x-y|^{4}}
 dxdy \\
& +  \iint  \frac{|v_{n,0}(x)-w^{J_1}_n(x)|^2 |v_{n,0}(y)-w^{J_1}_n(y)|^2}{|x-y|^{4}}
 dxdy -\sum^{J_1}_{j=1} \iint \frac{|V^l_j(-\frac{t_{j,n}}{\lambda^2_{j,n}},x)|^2 |V^l_j(-\frac{t_{j,n}}{\lambda^2_{j,n}},y)|^2}{|x-y|^{4}}
 dxdy \\
& + \iint  \frac{|w^J_n(x)-w^{J_1}_n(x)|^2 |w^J_n(y)-w^{J_1}_n(y)|^2}{|x-y|^{4}} \ dxdy -\iint  \frac{ |w^J_n(x)|^2 |w^J_n(y)|^2}{|x-y|^{4}} dxdy \\
& +  \sum^{J_1}_{j=J+1} \iint \frac{|V^l_j(-\frac{t_{j,n}}{\lambda^2_{j,n}},x)|^2 |V^l_j(-\frac{t_{j,n}}{\lambda^2_{j,n}},y)|^2}{|x-y|^{4}}
 dxdy -\iint  \frac{|w^J_n(x)-w^{J_1}_n(x)|^2 |w^J_n(y)-w^{J_1}_n(y)|^2}{|x-y|^{4}} \ dxdy   \Big|\\
 \leq & 3\epsilon,
\endaligned
\end{equation*}
this completes the proof.

\begin{lemma}\label{sc}
Let $\{ z_{0,n} \} \in \dot{H}^1$ be radial, with
\begin{equation*}
\big\|\nabla z_{0,n}\big\|_{L^2} < \big\| \nabla W \big\|_{L^2},
\quad E(z_{0,n}) \rightarrow E_c.
\end{equation*}
and with $\big\|e^{it\Delta}z_{0,n}\big\|_{X(\mathbb{R})}\geq
\delta>0$, where $\delta = \delta \big(\big\| \nabla W \big\|_{L^2}\big)$ is as in Proposition \ref{lwp}.
Let $V_{0,j}$ be as in Lemma \ref{pd}. Assume that one of the two
hypotheses holds
\begin{enumerate}
\item[$(1)$]
\begin{equation}\label{manydec}
\lim_{ \overline{n \rightarrow +\infty}}
E\big(V^l_1(-\frac{t_{1,n}}{\lambda^2_{1,n}})\big)<E_c.
\end{equation}
\item[$(2)$] After passing to a subsequence, we have that
\begin{equation}\label{onedec}
\lim_{ \overline{n \rightarrow +\infty}}
E\big(V^l_1(-\frac{t_{1,n}}{\lambda^2_{1,n}})\big)=E_c
\end{equation}
with $s_{1,n}=-\frac{t_{1,n}}{\lambda^2_{1,n}}\rightarrow s_* \in
[-\infty, \infty]$, and if $U_1$ is the nonlinear profile associated to
$\big(V_{0,1}, \{ s_{1,n} \} \big)$, we have that the maximal
interval of existence of $U_1$ is $I=(-\infty, +\infty)$ and $\big\|
U_1\big\|_{X(\mathbb{R})}<\infty$.
\end{enumerate}
Then, after passing to a subsequence, for $n$ large, if $z_n$ is the
solution of $(\ref{equ1})$ with data at $t=0$ equal to $z_{0,n}$, then
$(SC)(z_{0,n})$ holds.
\end{lemma}

{\bf Proof: Case 2 holds.}
Applying Lemma \ref{pd} to $\{ z_{0,n}\}$, we have
\begin{eqnarray}
 z_{0,n}(x)&
=& \sum^{J}_{j=1}\frac{1}{\lambda^{(d-2)/2}_{j,n}}
V^{l}_{j}\big(-\frac{t_{j,n}}{\lambda^2_{j,n}},
\frac{x}{\lambda_{j,n}} \big) + w^J_n \nonumber\\
\big\| \nabla W \big\|^2_{L^2} > \big\|\nabla z_{0,n }
\big\|^2_{L^2} & = & \sum^J_{j=1} \big\| \nabla V_{0,j}
\big\|^2_{L^2} + \big\|\nabla w^J_n \big\|^2_{L^2} +o_n(1), \nonumber \\
& = & \sum^J_{j=1} \big\| \nabla
V^l_j(-\frac{t_{j,n}}{\lambda^2_{j,n}})
\big\|^2_{L^2} + \big\|\nabla w^J_n \big\|^2_{L^2} +o_n(1),  \label{kesc}\\
E_c \leftarrow  E(z_{0,n}) & = & \sum^J_{j=1}
E(V^l_j(-\frac{t_{j,n}}{\lambda^2_{j,n}})) + E(w^J_n) +o_n(1).
\label{esc}
\end{eqnarray}

By $(\ref{kesc})$ and Corollary \ref{pe}, we have for every $1\leq j
\leq J$
\begin{equation*}
E(V^l_j(-\frac{t_{j,n}}{\lambda^2_{j,n}}))  \geq 0,  \quad
E(w^J_n)   \geq 0.
\end{equation*}
Using $(\ref{onedec})$ and $(\ref{esc})$, we have for every $2\leq j
\leq J$
\begin{equation*}
E(V^l_j(-\frac{t_{j,n}}{\lambda^2_{j,n}})) \rightarrow 0,  \quad
E(w^J_n)  \rightarrow 0, \quad \text{as}\ n\rightarrow +\infty.
\end{equation*}
Using Corollary \ref{cge}, we obtain that

\begin{equation*}
\aligned \sum^J_{j=2}\big\|\nabla V_{0,j}\big\|^2_{L^2} +
\big\|\nabla w^J_n\big\|^2_{L^2} =  \sum^J_{j=2}\big\|\nabla
V^l_j((-\frac{t_{j,n}}{\lambda^2_{j,n}}))\big\|^2_{L^2} +
\big\|\nabla w^J_n\big\|^2_{L^2} \rightarrow 0,\ \text{as}\ n \rightarrow +\infty.
\endaligned
\end{equation*}
Hence, we have for every $2\leq j \leq J$
\begin{equation*}
V_{0,j}\equiv 0, \quad \text{and} \ \ \big\|\nabla w^J_n\big\|_{L^2}
\rightarrow 0, \ \text{as}\ n\rightarrow +\infty.
\end{equation*}
Therefore,
\begin{equation*}
 z_{0,n}(x)
=\frac{1}{\lambda^{(d-2)/2}_{1,n}} V^{l}_{1}\big(s_{1,n},
\frac{x}{\lambda_{1,n}} \big) + w_n, \ \text{where}\ \big\| \nabla w_n
\big\|_{L^2} \rightarrow 0 \ \text{as} \ n\rightarrow +\infty.
\end{equation*}

Let $v_{0,n}=\lambda^{(d-2)/2}_{1,n} z_{0,n}(\lambda_{1,n}x)$,
$\widetilde{w}_n = \lambda^{(d-2)/2}_{1,n}w_n(\lambda_{1,n}x)$, we
have $\big\| \nabla v_{0,n} \big\|_{L^2} = \big\|\nabla
z_{0,n}\big\|_{L^2}<\big\|\nabla W\big\|_{L^2}$ and
\begin{equation*}
\aligned
 v_{0,n}(x) &=
V^{l}_{1}\big(s_{1,n}, x \big) + \widetilde{w}_n(x),  \ \text{where}\ \big\|
\nabla \widetilde{w}_n \big\|_{L^2} \rightarrow 0 \ \text{as} \ n\rightarrow +\infty.
\endaligned
\end{equation*}

Note that by the definition of nonlinear profile, we have
\begin{equation*}
\aligned \big\| \nabla U_1(s_{1,n})-\nabla
V^l_1(s_{1,n})\big\|_{L^2} \rightarrow 0 \ \text{as} \ n\rightarrow +\infty,
\endaligned
\end{equation*}
then
\begin{equation*}
\aligned
 v_{0,n}(x) &=
U_{1}\big(s_{1,n}, x \big) + \widetilde{\widetilde{w}}_n, \quad
\big\|
\nabla \widetilde{\widetilde{w}}_n \big\|_{L^2} \rightarrow 0. \\
E\big(U_1(s_{1,n})\big)&=E\big(V^l_1(s_{1,n})\big)+o_n(1) \rightarrow E_c, \\
\big\| \nabla U_1(s_{1,n})\big\|_{L^2}& = \big\|\nabla
V^l_1\big(s_{1,n}\big)\big\|_{L^2} + o_n(1) \\
& = \big\|\nabla
V_{0,1}\big\|_{L^2} + o_n(1) < \big\| \nabla W \big\|_{L^2}.
\endaligned
\end{equation*}

We now apply Proposition \ref{ltp} with $\widetilde{u}=U_1, e=0$ to
obtain that $(SC)(v_{0,n})$ holds, then this case follows from the
dilation invariance of (\ref{equ1}).

{\bf Case 1 holds.}
We first claim that
\begin{equation}\label{epl}
\aligned \lim_{ \overline{n \rightarrow +\infty}}
E\big(V^l_j(-\frac{t_{j,n}}{\lambda^2_{j,n}})\big) < E_c \quad \text{for}
\ \ j \geq 2.
\endaligned
\end{equation}

After passing to a subsequence, we assume that
\begin{equation}\label{ep1}
\lim_{n\rightarrow \infty}
E\big(V^l_1(-\frac{t_{1,n}}{\lambda^2_{1,n}})\big) < E_c.
\end{equation}
Applying Lemma \ref{pd} to $\{ z_{0,n}\}$, we have
\begin{eqnarray}
 z_{0,n}(x)&
=& \sum^{J}_{j=1}\frac{1}{\lambda^{\frac{d-2}{2}}_{j,n}}
V^{l}_{j}\big(-\frac{t_{j,n}}{\lambda^2_{j,n}},
\frac{x}{\lambda_{j,n}} \big) + w^J_n \nonumber\\
\big\| \nabla W \big\|^2_{L^2} > \big\|\nabla z_{0,n }
\big\|^2_{L^2} & = & \sum^J_{j=1} \big\| \nabla V_{0,j}
\big\|^2_{L^2} + \big\|\nabla w^J_n \big\|^2_{L^2} +o_n(1), \nonumber \\
& = & \sum^J_{j=1} \big\| \nabla
V^l_j(-\frac{t_{j,n}}{\lambda^2_{j,n}})
\big\|^2_{L^2} + \big\|\nabla w^J_n \big\|^2_{L^2} +o_n(1),  \label{kesc2}\\
E_c \leftarrow  E(z_{0,n}) & = & \sum^J_{j=1}
E(V^l_j(-\frac{t_{j,n}}{\lambda^2_{j,n}})) + E(w^J_n) +o_n(1).
\label{esc2}
\end{eqnarray}
By $(\ref{kesc2})$ and Corollary \ref{pe}, we have for every $1\leq j
\leq J$
\begin{equation*}
E(V^l_j(-\frac{t_{j,n}}{\lambda^2_{j,n}}))  \geq 0,  \quad
E(w^J_n)   \geq 0.
\end{equation*}

Note that
\begin{equation*}
\aligned \big\| \nabla
V^l_1(-\frac{t_{j,n}}{\lambda^2_{j,n}})\big\|^2_{L^2} & =  \big\|
\nabla V_{0,1}\big\|^2_{L^2} < \big\| \nabla W\big\|^2_{L^2},\\
E(V^l_1(-\frac{t_{j,n}}{\lambda^2_{j,n}})) & \leq E_c +o_n(1) < E(W).\\
\endaligned
\end{equation*}
Hence,  from Lemma \ref{slbv} and Lemma \ref{pd}, we have
\begin{equation*}
\aligned \int |\nabla V^l_1(-\frac{t_{j,n}}{\lambda^2_{j,n}}) |^2dx
& - \int \int
\frac{|V^l_1(-\frac{t_{j,n}}{\lambda^2_{j,n}},x)|^2|V^l_1(-\frac{t_{j,n}}{\lambda^2_{j,n}},y)|^2}{|x-y|^4}
dxdy
\geq \frac{\overline{\delta}}{2} \int |\nabla V^l_1(-\frac{t_{j,n}}{\lambda^2_{j,n}}) |^2dx, \\
E(V^l_1(-\frac{t_{j,n}}{\lambda^2_{j,n}}))  &= \frac14 \big\| \nabla
V^l_1(-\frac{t_{j,n}}{\lambda^2_{j,n}}) \big\|^2_{L^2} + \frac14
\big(  \int |\nabla V^l_1(-\frac{t_{j,n}}{\lambda^2_{j,n}}) |^2dx
\\
& \qquad \qquad  \qquad  \qquad \qquad  \qquad - \int \int
\frac{|V^l_1(-\frac{t_{j,n}}{\lambda^2_{j,n}},x)|^2|V^l_1(-\frac{t_{j,n}}{\lambda^2_{j,n}},y)|^2}{|x-y|^4}
dxdy \big)\\
& \geq C \big\| \nabla V^l_1(-\frac{t_{j,n}}{\lambda^2_{j,n}})
\big\|^2_{L^2}  = C \big\| \nabla
V_{0,1}\big\|^2_{L^2} \geq C
\alpha_0>0.
\endaligned
\end{equation*}

By (\ref{esc2}), we have
\begin{equation*}
\aligned E_c \leftarrow  E(z_{0,n}) \geq C \alpha_0 + \sum^J_{j=2}
E\big(V^l_j(-\frac{t_{j,n}}{\lambda^2_{j,n}})\big) + E(w^J_n) +o_n(1),
\endaligned
\end{equation*}
which implies the claim.

After passing to a subsequence, we can assume that for any $j\geq 1$
\begin{equation*}
\aligned \lim_{n\rightarrow \infty}
E\big(V^l_j((-\frac{t_{j,n}}{\lambda^2_{j,n}}))\big) \quad \text{exists}, \ \text{and}\quad
\lim_{n\rightarrow \infty}
-\frac{t_{j,n}}{\lambda^2_{j,n}}=\overline{s}_j \in [-\infty,
\infty].
\endaligned
\end{equation*}
If $U_j$ is the nonlinear profile associated to $\big(V_{0,j},
\{-\frac{t_{j,n}}{\lambda^2_{j,n}} \} \big)$, then by the definition
of nonlinear profile, for sufficiently large $n$, we obtain
\begin{equation*}
\aligned \big\|\nabla
U_j(-\frac{t_{j,n}}{\lambda^2_{j,n}})\big\|^2_{L^2} & = \big\|\nabla
V^{l}_j(-\frac{t_{j,n}}{\lambda^2_{j,n}})\big\|^2_{L^2} +o_n(1)<
\big\|\nabla
W\big\|^2_{L^2}, \\
E\big(U_j(-\frac{t_{j,n}}{\lambda^2_{j,n}})\big) &= E\big(V^l_j(-\frac{t_{j,n}}{\lambda^2_{j,n}})\big) + o_n(1)< E_c.
\endaligned
\end{equation*}
By the definition of $E_c$, we have that $U_j$ satisfies $(SC)$.
Moreover we also have $\big\|U_j\big\|_{W(\mathbb{R})}<
\infty$, and  we obtain from Corollary \ref{cge}
\begin{equation}\label{eqvke}
\aligned E\big(U_j(t)\big)\approx \big\|\nabla U_j(t)\big\|^2_{L^2} \approx
\big\|\nabla U_j(0)\big\|^2_{L^2},\ \forall\ t \in \mathbb{R}.
\endaligned
\end{equation}

On the other hand, we claim that there exists $j_0$ such that, for $j\geq j_0$
\begin{equation}\label{bpf}
\aligned \big\|U_j\big\|_{X(\mathbb{R})} \leq C \big\|\nabla
V_{0,j}\big\|_{L^2}.
\endaligned
\end{equation}
In fact, from (\ref{kesc2}), we have
\begin{equation*}
\aligned  \sum^J_{j=1} \big\| \nabla V_{0,j} \big\|^2_{L^2} \leq
\big\|\nabla z_{0,n } \big\|^2_{L^2} + o_n(1) \leq  \big\| \nabla W
\big\|^2_{L^2},
\endaligned
\end{equation*}
then there exists $j_0$, for $j\geq j_0$, such that $\big\|\nabla
V_{0,j}\big\|_{L^2} \leq \widetilde{\delta}$, where
$\widetilde{\delta}$ is so small that
$\big\|e^{it\Delta}V_{0,j}\big\|_{X(\mathbb{R})} \leq \delta$, with
$\delta$ as in Proposition \ref{lwp}. Note that
\begin{equation*}
\aligned U_j(t)=e^{it\Delta}V_{0,j}+i\int^t_{\overline{s}_j}
e^{i(t-s)\Delta}  \big(|x|^{-4}*|U_j|^2\big)(s,x)U_j(s,x) ds,
\endaligned
\end{equation*}
this together with the local wellposedness theory implies
\begin{equation*}
\aligned \big\|U_j\big\|_{X(\mathbb{R})} \leq C \big\|\nabla
V_{0,j}\big\|_{L^2}.
\endaligned
\end{equation*}

Since for sufficiently large $n$, we have
\begin{equation*}
\aligned z_{0,n}(x)& =&
\sum^{J(\epsilon_0)}_{j=1}\frac{1}{\lambda^{\frac{d-2}{2}}_{j,n}}
V^{l}_{j}\big(-\frac{t_{j,n}}{\lambda^2_{j,n}},
\frac{x}{\lambda_{j,n}} \big) + w^{J(\epsilon_0)}_n, \quad\big\|e^{it\Delta}w^{J(\epsilon_0)}_n\big\|_{{\cal Z}^1(\mathbb{R})}
\leq \epsilon_0.
\endaligned
\end{equation*}
Define the near-solution
\begin{equation*}
\aligned H_{n,\epsilon_0}(t,x)=\sum^{J(\epsilon_0)}_{j=1}
\frac{1}{\lambda^{\frac{d-2}{2}}_{j,n}}
U_j(\frac{t-t_{j,n}}{\lambda^2_{j,n}}, \frac{x}{\lambda_{j,n}}).
\endaligned
\end{equation*}
Then $H_{n,\epsilon_0}$ satisfies the following equation
\begin{equation*}
\aligned
   ( i\partial_t +  \Delta ) H_{n,\epsilon_0}(t, x)  & =  \sum^{J(\epsilon_0)}_{j=1}f\Big(\frac{1}{\lambda^{\frac{d-2}{2}}_{j,n}}
U_j(\frac{t-t_{j,n}}{\lambda^2_{j,n}}, \frac{x}{\lambda_{j,n}})\Big) \\
& =  f(H_{n,\epsilon_0}(t,x) ) + R_{n, \epsilon_0}(t,x)
\endaligned
\end{equation*}
where
\begin{equation*}
\aligned R_{n, \epsilon_0}(t,x)=
\sum^{J(\epsilon_0)}_{j=1}f\Big(\frac{1}{\lambda^{\frac{d-2}{2}}_{j,n}}
U_j(\frac{t-t_{j,n}}{\lambda^2_{j,n}}, \frac{x}{\lambda_{j,n}})\Big)-
f(H_{n,\epsilon_0}(t,x) ).
\endaligned
\end{equation*}
By the definition of the nonlinear profile $U_j$, we have
\begin{equation}\label{con4}
\aligned
 z_{0,n}(x) =&
\sum^{J(\epsilon_0)}_{j=1}\frac{1}{\lambda^{\frac{d-2}{2}}_{j,n}}
U_j\big(-\frac{t_{j,n}}{\lambda^2_{j,n}}, \frac{x}{\lambda_{j,n}}
\big) + \widetilde{w}^{J(\epsilon_0)}_n \nonumber\\
 = &H_{n, \epsilon_0}(0) + \widetilde{w}^{J(\epsilon_0)}_n, \quad
\big\|e^{it\Delta}\widetilde{w}^{J(\epsilon_0)}_n\big\|_{{\cal Z}^1(\mathbb{R})} \leq 2\epsilon_0  \ \text{for}\ n \gg 1.
\endaligned
\end{equation}
By the orthogonality property and (\ref{kesc2}), we have
\begin{equation}\label{con1}
\aligned
\big\|\nabla H_{n, \epsilon_0}(0)\big\|^2_{L^2} \leq  C \sum^{J(\epsilon_0)}_{j=1} \big\|\nabla
V^l_j(-\frac{t_{j,n}}{\lambda^2_{j,n}}) \big\|^2_{L^2} +o_n(1) \leq
 C  \big\|\nabla W\big\|^2_{L^2}.
 \endaligned
\end{equation}
In addition, we also have
\begin{equation*}
\aligned
\big\|H_{n, \epsilon_0}\big\|^6_{X(\mathbb{R})}  = & \int \Big\| \sum^{J(\epsilon_0)}_{j=1} \frac{1}{\lambda^{\frac{d-2}{2}}_{j,n}}
U_j(\frac{t-t_{j,n}}{\lambda^2_{j,n}}, \frac{x}{\lambda_{j,n}})
\Big\|^{6}_{L^{\frac{6d}{3d-8}}_x}dt \\
\leq & \int \Big(  \sum^{J(\epsilon_0)}_{j=1}  \Big\| \frac{1}{\lambda^{\frac{d-2}{2}}_{j,n}}
U_j(\frac{t-t_{j,n}}{\lambda^2_{j,n}}, \frac{x}{\lambda_{j,n}})
\Big\|_{L^{\frac{6d}{3d-8}}_x}\Big)^6 dt \\
 \leq &  \sum^{J(\epsilon_0)}_{j=1} \int    \Big\| \frac{1}{\lambda^{\frac{d-2}{2}}_{j,n}}
U_j(\frac{t-t_{j,n}}{\lambda^2_{j,n}}, \frac{x}{\lambda_{j,n}})
\Big\|^6_{L^{\frac{6d}{3d-8}}_x} dt \\
 & + C_{J(\epsilon_0)} \sum_{j\not =j'} \int \Big\| \frac{1}{\lambda^{\frac{d-2}{2}}_{j,n}}
U_j(\frac{t-t_{j,n}}{\lambda^2_{j,n}}, \frac{x}{\lambda_{j,n}})
\Big\|_{L^{\frac{6d}{3d-8}}_x}  \Big\| \frac{1}{\lambda^{\frac{d-2}{2}}_{j',n}}
U_{j'}(\frac{t-t_{j',n}}{\lambda^2_{j',n}}, \frac{x}{\lambda_{j',n}})
\Big\|^5_{L^{\frac{6d}{3d-8}}_x} dt \\
 =& I+II.
\endaligned
\end{equation*}
For the first term, from (\ref{kesc2}) and (\ref{bpf}), we have
\begin{equation*}
\aligned I & \leq \sum^{j_0}_{j=1}
\big\|U_{j}\big\|^{6}_{X(\mathbb{R})} +
\sum^{J(\epsilon_0)}_{j=j_0+1}
\big\|U_{j}\big\|^{6}_{X(\mathbb{R})}\\
& \leq \sum^{j_0}_{j=1} \big\|U_{j}\big\|^{6}_{X(\mathbb{R})} + C
\sum^{J(\epsilon_0)}_{j=j_0+1} \big\|\nabla V_{0,j}\big\|^6_{L^2}
\leq \sum^{j_0}_{j=1} \big\|U_{j}\big\|^{6}_{X(\mathbb{R})} + C
\big( \sum^{J(\epsilon_0)}_{j=j_0+1} \big\|\nabla V_{0,j}\big\|^2_{L^2}\big)^3 \\
& \leq \frac{C_0}{2},
\endaligned
\end{equation*}
where $C_0$ is independent of $J(\epsilon_0)$. For the second term,
we have from the orthogonality of $(\lambda_{j,n}, t_{j,n})$
\begin{equation*}
II \rightarrow 0 \quad \text{as} \ \ n\rightarrow \infty.
\end{equation*}
Hence, we obtain
\begin{equation}\label{con2}
\aligned
 \big\|H_{n, \epsilon_0}\big\|^6_{X(\mathbb{R})} \leq C_0, \ \text{for} \ n \ \text{sufficiently large},
\endaligned
\end{equation}
where $C_0$ is independent of $J(\epsilon_0)$.

Note that $\big\|U_j\big\|_{X(\mathbb{R})}<\infty$ and $\big\|
U_j\big\|_{W(\mathbb{R})}<\infty$, using the orthogonality of $(\lambda_{j,n}, t_{j,n})$ again, we have that
\begin{equation}\label{con3}
\aligned
\big\|R_{n,\epsilon_0}(t,x)\big\|_{L^{\frac32}(\dot{H}^{1,\frac{6d}{3d+4}})}
\rightarrow 0 \quad \text{as} \ \ n\rightarrow \infty.
\endaligned
\end{equation}
Last, for sufficiently large $n$, we have
\begin{equation}\label{con5}
\aligned
\big\|\nabla z_{0,n}-\nabla H_{n,\epsilon_0}(0)\big\|_{L^2} \leq &
\big\|\widetilde{w}^{J(\epsilon_0)}_n\big\|_{L^2} \leq  \big\|w^{J(\epsilon_0)}_n\big\|_{L^2} +o_n(1) \leq \big\|\nabla W\big\|^2_{L^2}.
\endaligned
\end{equation}
Combining Proposition \ref{ltp}, Remark \ref{stab1} with (\ref{con1})-(\ref{con5}), we
obtain that $(SC)(z_{0,n})$ holds.

\begin{proposition}[Existence of a critical solution]\label{ecs}
There exists a radial solution $u_c$ of $(\ref{equ1})$ in $\dot{H}^1$
with data $u_{c,0}$ and maximal interval of existence $I$  such that
\begin{equation*}
\big\|\nabla u_{c,0} \big\|_{L^2} < \big\|\nabla W \big\|_{L^2}, \quad E(u_{c,0}) =
E_c
\end{equation*}
and
\begin{equation*}
\big\|u_c \big\|_{X(I)}=+\infty.
\end{equation*}
\end{proposition}

{\bf Proof: } By the definition of $E_c$ and the assumption that
$E_c < E(W)$, we can find $u_{0,n} \in \dot{H}^1$ radial, with
$\big\| \nabla u_{0,n}\big\|_{L^2} < \big\|\nabla W \big\|_{L^2}$,
$E(u_{0,n}) \searrow E_c$, and such that if $u_n$ is the solution of
(\ref{equ1}) with data $u_{0,n}$ at $t=0$ and maximal interval of
existence $I_n=(-T_{-}(u_{0,n}), T_{+}(u_{0,n}))$, then
\begin{equation*}
\aligned \big\|e^{it\Delta} u_{0,n}\big\|_{X(\mathbb{R})} \geq
\delta \quad \text{as Proposition }  \ref{lwp},\ \text{and}\
\big\|u_n\big\|_{S(I_n)}=+\infty.
\endaligned
\end{equation*}

Note that $E_c < E(W)$, then there exists $\delta_0>0$, so that for
sufficiently large $n$, we have $E(u_{0,n}) < (1-\delta_0) E(W)$. By Proposition
\ref{dlbv}, we can find $\overline{\delta}$ so that
\begin{equation*}
\aligned \big\|\nabla u_n(t)\big\|^2_{L^2} \leq
(1-\overline{\delta}) \big\|\nabla W\big\|^2_{L^2}, \forall \ t\in
I_n.
\endaligned
\end{equation*}

Applying Lemma \ref{pd} to $\{ u_{0,n} \}$, we have
\begin{equation*}
\aligned
 u_{0,n} =&\sum^{J}_{j=1}\frac{1}{\lambda^{(d-2)/2}_{j,n}}
V^{l}_{j}\big(-\frac{t_{j,n}}{\lambda^2_{j,n}},
\frac{x}{\lambda_{j,n}} \big) + w^J_n \nonumber\\
\big\|V_{0,1}\big\|_{\dot{H}^1}  \geq & \alpha_0(A)>0, \quad
\lim_{J \rightarrow \infty} \big[ \lim_{n \rightarrow \infty} \big\|
e^{it\Delta} w^J_n \big\|_{X(\mathbb{R})} \big]=0,
\endaligned
\end{equation*}
\begin{eqnarray}
(1-\overline{\delta}) \big\|\nabla W\big\|^2_{L^2} \geq \big\|\nabla u_{0,n } \big\|^2_{L^2} & =& \sum^J_{j=1} \big\|\nabla
V_{0,j}
\big\|^2_{L^2} + \big\|\nabla w^J_n \big\|^2_{L^2} +o_n(1), \nonumber\\
&= & \sum^J_{j=1} \big\|\nabla
V^l_j(-\frac{t_{j,n}}{\lambda^2_{j,n}})
\big\|^2_{L^2} + \big\|\nabla w^J_n \big\|^2_{L^2} +o_n(1), \label{kiecs}\\
E_c \swarrow E(u_{0,n}) & = & \sum^J_{j=1}
E(V^l_j(-\frac{t_{j,n}}{\lambda^2_{j,n}})) + E(w^J_n) +o_n(1).
\label{encs}
\end{eqnarray}

Because of (\ref{kiecs}), we have that
\begin{equation*}
\aligned
\big\|\nabla V^l_j(-\frac{t_{j,n}}{\lambda^2_{j,n}}) \big\|^2_{L^2}
\leq (1-\frac{\overline{\delta}}{2}) \big\|\nabla W\big\|^2_{L^2},
\quad \big\|\nabla w^J_n \big\|^2_{L^2} \leq
(1-\frac{\overline{\delta}}{2}) \big\|\nabla W\big\|^2_{L^2}, \quad \text{for}\ n\gg 1.
\endaligned
\end{equation*}
From Corollary \ref{pe}, it follows that
\begin{equation*}
E(V^l_j(-\frac{t_{j,n}}{\lambda^2_{j,n}})) \geq 0,\quad E(w^J_n) \geq 0.
\end{equation*}
By (\ref{encs}), we have that
\begin{equation*}
E(V^l_1(-\frac{t_{1,n}}{\lambda^2_{1,n}})) \leq E(u_{0,n})+o_n(1),
\end{equation*}
therefore,
\begin{equation*}
\liminf_{n\rightarrow
\infty}E(V^l_1(-\frac{t_{1,n}}{\lambda^2_{1,n}})) \leq E_c.
\end{equation*}

Note that $(SC)(u_{0,n})$ does not hold, we have from Lemma \ref{sc}
\begin{equation*}
\aligned
\liminf_{ n \rightarrow
\infty} E(V^l_1(-\frac{t_{1,n}}{\lambda^2_{1 ,n}})) = E_c.
\endaligned
\end{equation*}

Arguing as in the proof of Case 2,  Lemma \ref{sc}, we see that $\displaystyle
\liminf_{n\rightarrow
\infty}E(V^l_1(-\frac{t_{1,n}}{\lambda^2_{1 ,n}})) = E_c $ and $E_c
< E(W)$ imply that $J=1$ and $\big\| \nabla w^J_n\big\|_{L^2} \rightarrow
0$ as $n\rightarrow +\infty$.

Thus
\begin{equation*}
u_{0,n}(x) =\frac{1}{\lambda^{(d-2)/2}_{1,n}}
V^{l}_{1}\big(-\frac{t_{1,n}}{\lambda^2_{1,n}},
\frac{x}{\lambda_{1,n}} \big) + w_n(x),\quad \big\|\nabla w_n
\big\|_{L^2} \rightarrow 0 \ \text{as}\ n\rightarrow +\infty.
\end{equation*}

Let
\begin{equation*}
v_{0,n}(x) = \lambda^{(d-2)/2}_{1,n} u_{0,n}(\lambda_{1,n}x), \quad
\widetilde{w}_n(x)= \lambda^{(d-2)/2}_{1,n} w_n(\lambda_{1,n}x),
\end{equation*}
then
\begin{equation*}
v_{0,n}(x) = V^{l}_{1}\big(-\frac{t_{1,n}}{\lambda^2_{1,n}}, x \big)
+ \widetilde{w}_n(x), \quad \big\|\nabla \widetilde{w}_n
\big\|_{L^2} \rightarrow 0 \ \text{as}\ n \rightarrow +\infty.
\end{equation*}

Let $U_1$ be the nonlinear profile associated to $(V_{0,1},
-\frac{t_{1,n}}{\lambda^2_{1,n}})$ and let $I_1$ be its maximal
interval of existence. By the definition of the nonlinear profile,
we have for $s_n=-\frac{t_{1,n}}{\lambda^2_{1,n}}$
\begin{equation*}
\aligned \big\|\nabla U_1(s_n)\big\|^2_{L^2} & = \big\|\nabla
V^l_1(s_n)\big\|^2_{L^2}+o_n(1) < \big\|\nabla W\big\|^2_{L^2},\\
E(U_1(s_n)) & = E(V^l_1(s_n))+o_n(1) = E_c +o_n(1).
\endaligned
\end{equation*}
Let's fix $s_{*} \in I_1$. then from the conservation of energy and
Proposition \ref{dlbv}, we have
\begin{equation*}
\aligned \big\|\nabla U_1(s_{*})\big\|^2_{L^2}  <  \big\|\nabla
W\big\|^2_{L^2}, \quad  E(U_1(s_{*}))  =  E_c.
\endaligned
\end{equation*}

If $\big\|U_1 \big\|_{X(I_1)} < +\infty$, Proposition \ref{buc}
implies that $I_1 = (-\infty, +\infty)$, then $(SC)(u_{0,n})$ holds
from Lemma \ref{sc}, this obtains a contradiction. Thus
\begin{equation*}
\big\|U_1 \big\|_{X(I_1)} = +\infty.
\end{equation*}
This completes the proof.

\begin{proposition}[Pre-compactness of the flow of the critical
solution]\label{ccs} Let $u_c$ be as in Proposition \ref{ecs}, and
that $\big\|u_c\big\|_{X(I_{+})}=+\infty$, where $I_{+}=(0,
+\infty)\cap I$. Then for $ t\in I_{+}$, there exists $\lambda(t)\in
\mathbb{R}^{+}$, such that $K$ is precompact in $\dot{H}^1$ where
\begin{equation*}
K=\Big\{ v(t,x), v(t,x)=\frac{1}{\lambda(t)^{\frac{d-2}{2}}}u_c(t,
\frac{x}{\lambda(t)}), t\in I_{+}  \Big\}.
\end{equation*}
\end{proposition}

{\bf Proof: } For brevity of notation,
let us set $u(t,x)=u_c(t,x)$.  We argue by contradiction. If not, there exist $\eta_0>0$ and a
sequence $\{t_n\}^{\infty}_{n=1}, t_n\geq 0$ such that, for all
$\lambda_0 \in \mathbb{R}^{+}$, we have
\begin{equation}\label{contra1}
\aligned \Big\|\frac{1}{\lambda^{(d-2)/2}_0}u(t_n,
\frac{x}{\lambda_0}) - u(t_{n'},x) \Big\|_{\dot{H}^1} \geq
\eta_0,\quad \text{for}\ \ n\not= n'.
\endaligned
\end{equation}

After passing to a subsequence, we assume that $t_n \rightarrow
\overline{t}\in [0, T_{+}(u_0)]$. By taking $\lambda_0=1$ in
(\ref{contra1}) and the continuity of the flow $u(t)$ in $\dot{H}^1$, we
must have
\begin{equation*}
\aligned \overline{t}= T_{+}(u_0).
\endaligned
\end{equation*}
In addition, from Proposition \ref{lwp}, we also have
\begin{equation}\label{contra2}
\aligned \big\| e^{it\Delta}u(t_n) \big\|_{S(0,+\infty)} \geq \delta.
\endaligned
\end{equation}

Applying Lemma \ref{pd} to $v_{0,n}=u(t_n)$, we have
\begin{equation*}
u(t_n,x)=\sum^{J}_{j=1}\frac{1}{\lambda^{(d-2)/2}_{j,n}}
V^{l}_{j}\big(-\frac{t_{j,n}}{\lambda^2_{j,n}},
\frac{x}{\lambda_{j,n}} \big) + w^J_n(x),
\end{equation*}
with
\begin{equation*}
\aligned \big\|\nabla u(t_n) \big\|^2_{L^2} & =\sum^J_{j=1}
\big\|\nabla V_{0,j}
\big\|^2_{L^2} + \big\|\nabla w^J_n \big\|^2_{L^2} +o_n(1), \\
E(u(t_n)) & = \sum^J_{j=1}
E(V^l_j(-\frac{t_{j,n}}{\lambda^2_{j,n}})) + E(w^J_n) +o_n(1).
\endaligned
\end{equation*}

Arguing as in the proof of Proposition \ref{ecs}, we see that
\begin{equation*}
\aligned
\liminf_{ n \rightarrow
\infty} E(V^l_1(-\frac{t_{1,n}}{\lambda^2_{1 ,n}})) = E_c,
\endaligned
\end{equation*}
this implies that $J=1$, i. e.
\begin{equation}\label{profile2}
\aligned u(t_n)= \frac{1}{\lambda^{(d-2)/2}_{1,n}}
V^{l}_{1}\big(-\frac{t_{1,n}}{\lambda^2_{1,n}},
\frac{x}{\lambda_{1,n}} \big) + w_n, \  \lim_{n\rightarrow +\infty}\big\| w_n\big\|_{\dot{H}^1} \rightarrow 0 .
\endaligned
\end{equation}

The next step is to show that
\begin{equation*}
\aligned s_n=-\frac{t_{1,n}}{\lambda^2_{1,n}} \ \ \text{must be
bounded}.
\endaligned
\end{equation*}
Notice that we have
\begin{equation*}
\aligned e^{it\Delta}u(t_n)= \frac{1}{\lambda^{(d-2)/2}_{1,n}}
V^{l}_{1}\big(\frac{t- t_{1,n}}{\lambda^2_{1,n}},
\frac{x}{\lambda_{1,n}} \big) + e^{it\Delta} w_n,
\endaligned
\end{equation*}
with $ \big\| e^{it\Delta} w_n\big\|_{X(\mathbb{R})} <
\frac{\delta}{2} $ for $n$ sufficiently large.

Assume that $\frac{t_{1,n}}{\lambda^2_{1,n}}\leq -C_0$ for $n$
large, $C_0$ a large positive constant. Since
\begin{equation*}
\aligned \Big\| \frac{1}{\lambda^{(d-2)/2}_{1,n}}
V^{l}_{1}\big(\frac{t- t_{1,n}}{\lambda^2_{1,n}},
\frac{x}{\lambda_{1,n}} \big)  \Big\|_{X(0, +\infty)} \leq
\big\|V^l_1 \big\|_{X(C_0, \infty)} < \frac{\delta}{2}
\endaligned
\end{equation*}
for $C_0$ large, we get for $n$ large
\begin{equation*}
\aligned \big\| e^{it\Delta}u(t_n) \big\|_{X(0,+\infty)} < \delta,
\endaligned
\end{equation*}
which is a contradiction to (\ref{contra2}).

On the other hand, if $\frac{t_{1,n}}{\lambda^2_{1,n}}\geq C_0$ for
$n$ large, we have
\begin{equation*}
\aligned \Big\| \frac{1}{\lambda^{(d-2)/2}_{1,n}}
V^{l}_{1}\big(\frac{t- t_{1,n}}{\lambda^2_{1,n}},
\frac{x}{\lambda_{1,n}} \big)  \Big\|_{X(-\infty, 0)} \leq
\big\|V^l_1 \big\|_{X(-\infty, -C_0)} < \frac{\delta}{2}
\endaligned
\end{equation*}
for $C_0$ large. Hence,
\begin{equation*}
\aligned \big\| e^{it\Delta}u(t_n) \big\|_{X(-\infty, t_n)} \leq
\delta
\endaligned
\end{equation*}
for $n$ large, Proposition \ref{lwp} now gives
\begin{equation*}
\aligned \big\|u\big\|_{X(-\infty, t_n)} \leq 2 \delta.
\endaligned
\end{equation*}
Since $t_n \rightarrow \overline{t}= T_{+}(u_0)$, we also obtain a
contradiction.

Hence
\begin{equation*}
\aligned \Big|-\frac{t_{1,n}}{\lambda^2_{1,n}} \Big| \leq C_0,
\endaligned
\end{equation*}
after passing to a subsequence, we can assume that
\begin{equation*}
\aligned \frac{t_{1,n}}{\lambda^2_{1,n}} \rightarrow t_0 \in
(-\infty, +\infty).
\endaligned
\end{equation*}

On the other hand, by (\ref{contra1}) and (\ref{profile2}), we obtain that for $n
\not= n'$ large,
\begin{equation*}
\aligned \Big\| \frac{1}{\lambda^{(d-2)/2}_0}
\frac{1}{\lambda^{(d-2)/2}_{1,n}}
V^l_1\big(-\frac{t_{1,n}}{\lambda^2_{1,n}},
\frac{\frac{x}{\lambda_0}}{\lambda_{1,n}} \big) -
  \frac{1}{\lambda^{(d-2)/2}_{1,n'}}
V^l_1\big(-\frac{t_{1,n'}}{\lambda^2_{1,n'}},
\frac{x}{\lambda_{1,n'}} \big)\Big\|_{\dot{H}^1} \geq
\frac{\eta_0}{2},
\endaligned
\end{equation*}
or
\begin{equation*}
\aligned \Big\|
\big(\frac{\lambda_{1,n'}}{\lambda_0\lambda_{1,n}}\big)^{(d-2)/2}
V^l_1\big(-\frac{t_{1,n}}{\lambda^2_{1,n}},
\frac{\lambda_{1,n'}}{\lambda_0\lambda_{1,n}}y \big) -
V^l_1\big(-\frac{t_{1,n'}}{\lambda^2_{1,n'}}, y
\big)\Big\|_{\dot{H}^1} \geq \frac{\eta_0}{2}.
\endaligned
\end{equation*}

Letting
\begin{equation*}
\aligned \lambda_0=\frac{\lambda_{1,n'}}{\lambda_{1,n}},
\endaligned
\end{equation*}
we will obtain a contradiction because of the continuity of the linear flow $V^l_1(t,x)$ in $\dot{H}^1$ and
\begin{equation*}
\aligned -\frac{t_{1,n}}{\lambda^2_{1,n}} \rightarrow t_0, \quad
\text{and} \ \ -\frac{t_{1,n'}}{\lambda^2_{1,n'}} \rightarrow t_0.
\endaligned
\end{equation*}
This completes the proof.

%%%%%%%%%%%%%%%%%%%%%%%%%%%%%%%%%%%%%%%%%%%%%%%%%%%%%%%%%%%%%%%%%%%%%%%%%%%%%%%%%%%
\section{Rigidity theorem}
\setcounter{section}{5} \setcounter{equation}{0} In this section, we
will prove main theorem.

\begin{theorem}\label{rigthm}
Assume that $u_0 \in \dot{H}^1$ is radial and satisfies
\begin{equation*}
\aligned E(u_0)<E(W), \quad \big\|\nabla u_0\big\|_{L^2} <
\big\|\nabla W\big\|_{L^2}.
\endaligned
\end{equation*}
Let $u$ be the solution of $(\ref{equ1})$ with maximal interval of
existence $(-T_{-}(u_0), T_{+}(u_0))$. Assume that there exists
$\lambda(t)>0$, for $t\in [0, T_{+}(u_0))$, with the property that
\begin{equation*}
\aligned K=\big\{v(t,x)=\frac{1}{\lambda(t)^{\frac{d-2}{2}}} u(t,
\frac{x}{\lambda(t)}), t\in [0, T_{+}(u_0)) \big\}
\endaligned
\end{equation*}
is precompact in $\dot{H}^1$. Then
$T_{+}(u_0)=+\infty, u_0\equiv 0.$
\end{theorem}

We start out with a special case of the strengthened form of Theorem
\ref{rigthm}
\begin{proposition}\label{strengthm}
Assume that $u, v, \lambda(t)$ are as in Theorem \ref{rigthm}, and
that $\lambda(t) \geq A_0 >0$. Then the conclusion of Theorem
\ref{rigthm} holds.
\end{proposition}

First we collect some useful facts:

\begin{lemma}\label{usefulfacts}
Let $u,v$ be as in Theorem \ref{rigthm}.
\begin{enumerate}
\item[$(1)$] Let $\delta_0 > 0$ be such that $E(u_0) \leq
(1-\delta_0)E(W)$. Then there exists $\overline{\delta}>0$ such that for all $t\in [0, T_{+}(u_0))$, we have
\begin{equation}\label{use1}
\aligned \int \big| \nabla u(t)\big|^2 dx & \leq
(1-\overline{\delta})
 \int \big| \nabla W \big|^2dx, \\
\int |\nabla u(t,x) |^2dx  & -  \int  \int
\frac{|u(t,x)|^2|u(t,y)|^2}{|x-y|^4}
dxdy \geq \frac{\overline{\delta}}{2} \int |\nabla u |^2dx,\\
 \int \big| \nabla u(t)\big|^2 dx  \approx  E&( u(t)) = E(u_0)  \approx  \int \big| \nabla u_0 \big|^2 dx.
\endaligned
\end{equation}
\item[$(2)$] For all $t\in [0, T_{+}(u_0))$, we have
 \begin{equation*}
\aligned \big\| v(t,x)\big\|^2_{L^{2^*}}  \leq C_1 \int \big|\nabla v(t,x) \big|^2 dx  & \leq C_2 \int \big|
\nabla W(x)\big|^2dx.
\endaligned
\end{equation*}
\item[$(3)$] For each $\epsilon$, there exists $R(\epsilon)>0$,
such that for $ t\in [0, T_{+}(u_0))$,  we have
\begin{equation} \label{use3}
\aligned \int_{|x| > R(\epsilon)} \big| \nabla v(t,x) \big|^2 +
\big| v(t,x) \big|^{2^*} + \frac{\big| v(t,x)\big|^2}{|x|^2}dx +  \iint_{\Omega}
\frac{|v(t,x)|^2|v(t,y)|^2}{|x-y|^4}dxdy  \leq
\epsilon,
\endaligned
\end{equation}
where
\begin{equation*}
\aligned
\Omega=\big\{(x,y)\in \mathbb{R}^d\times \mathbb{R}^d; |x|>
R(\epsilon) \big\} \cup \big\{(x,y)\in \mathbb{R}^d\times
\mathbb{R}^d; |y|> R(\epsilon) \big\}.
\endaligned
\end{equation*}
\end{enumerate}
\end{lemma}
{\bf Proof: } From the property of $K$, we can easily verify them.

{\bf Proof of Proposition \ref{strengthm}: } We split the proof into
two cases, the finite time blowup for $u$ and the infinite time
of existence for $u$.

{\bf Case 1: $T_{+}(u_0)< + \infty$. } We claim that
\begin{equation*}
\aligned \lambda(t) \rightarrow \infty \quad \text{as}\ \
t\rightarrow T_{+}(u_0).
\endaligned
\end{equation*}
Its proof is analogue to the proof of Proposition $5.3$ in
\cite{KeM06} and Corollary $1.10$ in \cite{Ker06}. If not, there
exists $t_i \nearrow T_{+}(u_0)$ with $\lambda(t_i) \rightarrow
\lambda_0 \in [A_0, + \infty)$.

Let
\begin{equation*}
\aligned v_i(x)= \frac{1}{\lambda(t_i)^{\frac{d-2}{2}}}u\big(t_i,
\frac{x}{\lambda(t_i)} \big)
\endaligned
\end{equation*}
from the compactness of $\overline{K}$, there exists $v(x)\in
\dot{H}^1$ with
\begin{equation*}
\aligned v_i \rightarrow v \quad \text{in} \ \ \dot{H}^1,
\endaligned
\end{equation*}
Thus, we have
\begin{equation*}
\aligned u\big(t_i, x \big) = \lambda(t_i)^{\frac{d-2}{2}}
v_i(\lambda(t_i) x ) \rightarrow \lambda_0^{\frac{d-2}{2}}
v(\lambda_0 x) \quad \text{in} \ \ \dot{H}^1,
\endaligned
\end{equation*}
Let $h(t,x)$ be the solution of (\ref{equ1}) with data
$\lambda_0^{\frac{d-2}{2}} v(\lambda_0 x)$ at time $T_{+}(u_0)$ in
an interval $\big(T_{+}(u_0)-\delta, T_{+}(u_0)+\delta\big)$ with
\begin{equation*}
\aligned \big\|h\big\|_{X\big( (T_{+}(u_0)-\delta,
T_{+}(u_0)+\delta) \big)} < \infty.
\endaligned
\end{equation*}
Let $h_i(t,x)$ be the solution with data at $T_{+}(u_0)$ equal to
$u(t_i, x)$. Then the local well-posedness theory and Remark \ref{ct} guarantee that
\begin{equation*}
\aligned \sup_{i} \big\| h_i(t,x) \big\|_{X\big(
(T_{+}(u_0)-\frac{\delta}{2}, T_{+}(u_0)+\frac{\delta}{2}) \big)} <
\infty.
\endaligned
\end{equation*}
Since $h_i(t,x)= u(t+t_i-T_{+}(u_0),x )$ and $T_{+}(u_0)<\infty$, It
gives a contradiction with Proposition \ref{buc}.

Now let $\varphi \in C^{\infty}_0(\mathbb{R}^d)$ be radial, and
\begin{equation*}
\aligned \varphi(x)= \left\{\begin{array}{rl} 1, & \text{for} \ \ |x|\leq 1; \\
0, & \text{for} \ \ |x|\geq 2 \end{array}\right.
\endaligned
\end{equation*}
Set
\begin{equation*}
\aligned \varphi_R(x)=\varphi(\frac{x}{R}).
\endaligned
\end{equation*}

Define
\begin{equation*}
\aligned y_R(t) = \int \varphi_R(x) \big| u(t,x) \big|^2dx, \quad t
\in [0, T_{+}(u_0)).
\endaligned
\end{equation*}
From Lemma \ref{usefulfacts} and Lemma \ref{virialidentity}, we then
have
\begin{equation}\label{iden1}
\aligned \big| y'_R(t) \big| & \lesssim \int \Big| u(t,x) \nabla
u(t,x) \nabla\big(\varphi_R(x) \big) \Big|dx \\
& \lesssim  \big\|\nabla u(t) \big\|_{L^2} \Big\|
\frac{u(t,x)}{|x|}\Big\|_{L^2} \lesssim \big\| \nabla
W(x)\big\|^2_{L^2}.
\endaligned
\end{equation}

On the other hand, we also have
\begin{equation}\label{iden2}
\aligned \forall \ R>0, \quad \int_{|x|<R} \big|u(t,x) \big|^2dx
\rightarrow 0 \quad \text{as}\ \ t \rightarrow T_{+}(u_0).
\endaligned
\end{equation}
Indeed, since $u(t,x) = \lambda(t)^{\frac{d-2}{2}} v(t,
\lambda(t)x)$, we have from H\"{o}lder's inequality
\begin{equation*}
\aligned \int_{|x|<R} \big|u(t,x) \big|^2dx & = \lambda(t)^{-2}
\int_{|y|<\lambda(t)R}  \big|v(t,y) \big|^2dy \\
& = \lambda(t)^{-2} \int_{|y|<\epsilon \lambda(t)R}  \big|v(t,y)
\big|^2dy + \lambda(t)^{-2} \int_{\epsilon \lambda(t)R \leq  |y|
\leq \lambda(t)R} \big|v(t,y) \big|^2dy \\
& \leq \lambda(t)^{-2} \big( \epsilon \lambda(t)R \big)^{2} \big\|v(t,x)\big\|^2_{L^{2^*}} + \lambda(t)^{-2} \big(
\lambda(t)R \big)^{2}
\big\|v(t,x)\big\|^2_{L^{2^*}(|x|\geq \epsilon \lambda(t)R )} \\
& = C_3 \big( \epsilon R \big)^2 \int \big| \nabla W(x)\big|^2dx +
R^2 \big\|v(t,x)\big\|^2_{L^{2^*}(|x|\geq \epsilon \lambda(t)R )}.
\endaligned
\end{equation*}
The first term is small with $\epsilon$. Lemma \ref{usefulfacts}
implies that the second term tends to $0$ as $t$ tends to
$T_{+}(u_0)$.

From (\ref{iden2}), we have
\begin{equation}\label{iden3}
\aligned y_R(t) \rightarrow 0 \quad \text{as} \ \ t \rightarrow
T_{+}(u_0).
\endaligned
\end{equation}

From (\ref{iden1}) and (\ref{iden3}), we have
\begin{equation*}
\aligned y_R(0) & \leq y_R(T_{+}(u_0))+ C \ T_{+}(u_0)
\int\big|\nabla W(x)\big|^2dx\\
& = C \ T_{+}(u_0) \int\big|\nabla W(x)\big|^2dx
\endaligned
\end{equation*}
where $y_R(T_{+}(u_0))$ denotes  $\displaystyle
\lim_{t\nearrow T_{+}(u_0)} y_R(t)$.

Thus, letting $R \rightarrow +\infty$, we obtain
\begin{equation*}
\aligned u_0 \in L^2(\mathbb{R}^d).
\endaligned
\end{equation*}
Arguing as before, we have
\begin{equation*}
\aligned \big| y_R(t)\big|= \big| y_R(t) - y_R(T_{+}(u_0))\big| &
\leq C \ \big( T_{+}(u_0) -t\big) \int\big|\nabla W(x)\big|^2dx.
\endaligned
\end{equation*}
Letting $R \rightarrow +\infty$, we have
\begin{equation*}
\aligned \big\| u(t)\big\|^2_{L^2} \leq C \ \big( T_{+}(u_0)
-t\big) \int\big|\nabla W(x)\big|^2dx.
\endaligned
\end{equation*}
By the conservation of mass, this implies
\begin{equation*}
\aligned u_0 \equiv 0
\endaligned
\end{equation*}
which is in contradiction with $T_{+}(u_0)< +\infty$.

 {\bf Case 2: $T_{+}(u_0)=+ \infty$. }
On one hand, from  $u(t,x)=\lambda(t)^{\frac{d-2}{2}} v(t, \lambda(t) x)$ and Lemma \ref{usefulfacts}, we have for each $\epsilon>0$, there
exists $R(\epsilon)>0$ such that
\begin{equation}\label{decay1}
\aligned \int_{|x|> R(\epsilon)} \frac{|u(t,x)|^2}{|x|^2} dx +
\int_{|x|> R(\epsilon)} |\nabla u(t,x)|^2 dx +  \iint_{\Omega}
\frac{|u(t,x)|^2|u(t,y)|^2}{|x-y|^4}dxdy \leq \epsilon,
\endaligned
\end{equation}
where
\begin{equation*}
\Omega=\big\{(x,y)\in \mathbb{R}^d\times \mathbb{R}^d; |x|>
R(\epsilon) \big\} \cup \big\{(x,y)\in \mathbb{R}^d\times
\mathbb{R}^d; |y|> R(\epsilon) \big\}.
\end{equation*}

On the other hand, from Lemma \ref{usefulfacts}, and (\ref{decay1}), there exists $R$ such that, for all $t\in [0,
+\infty)$
\begin{equation}\label{decay2}
\aligned 8\int_{|x|\leq R} |\nabla u(t,x)|^2 dx -8 \iint_{\Omega_1}
\frac{|u(t,x)|^2|u(t,y)|^2}{|x-y|^4}dxdy \geq C_{\delta_0} \int |
\nabla u_0(x)|^2dx,
\endaligned
\end{equation}
where
\begin{equation*}
\aligned \Omega_1 & =\big\{ (x,y)\in \mathbb{R}^d\times
\mathbb{R}^d; |x|\leq R, |y|\leq R \big\}.
\endaligned
\end{equation*}

Now let $\varphi \in C^{\infty}_0(\mathbb{R}^d)$ be radial, and
\begin{equation*}
\aligned \varphi(x)= \left\{\begin{array}{rl} |x|^2, & \text{for} \ \ |x|\leq 1; \\
0, & \text{for} \ \ |x|\geq 2, \end{array}\right.
\endaligned
\end{equation*}
Set
\begin{equation*}
\aligned \varphi_R(x)=R^2\varphi(\frac{x}{R}).
\endaligned
\end{equation*}

Define
\begin{equation*}
\aligned z_R(t) = \int \varphi_R(x) \big| u(t,x) \big|^2dx, \quad t
\in [0, T_{+}(u_0)).
\endaligned
\end{equation*}
We then have
\begin{equation}\label{con6}
\aligned \big| z'_R(t) \big|  \leq C  R^2 \int
\big| \nabla u_0 \big|^2dx, & \quad \text{for}\ t > 0 , \\
 z''_R(t) \geq
C_{ \delta_0} \int \big| \nabla u_0 \big|^2dx, \quad &\quad  \text{for} \  R \ \text{large enough}, t > 0.
\endaligned
\end{equation}
In fact, from Lemma \ref{usefulfacts} and Lemma
\ref{virialidentity}, we have
\begin{equation*}
\aligned \big| z'_R(t) \big| & \leq 2 R  \int \Big|
\overline{u}(t,x) \nabla u(t,x) \nabla \varphi(\frac{x}{R}) \Big| dx
\\
& \leq C  R \int_{|x| \leq 2R} \big|u\big| \big| \nabla u\big|
 dx  \leq C  R^2 \big\| \nabla u(t,x) \big\|_{L^2}
\Big\|\frac{|u|}{|x|} \Big\|_{L^2} \leq C R^2 \int \big| \nabla
u_0 \big|^2dx.
\endaligned
\end{equation*}
On the other hand, from Lemma \ref{virialidentity}, (\ref{decay1})
and (\ref{decay2}), we have for sufficiently large $R$
\begin{equation*}
\aligned z''_R(t) & = -\int\triangle\triangle \varphi
(\frac{x}{R})\frac{|u|^{2}}{R^2}dx+4\mathrm{Re}\int
\varphi_{jk}\overline{u}_{j}u_{k}dx\\
& \quad -4\mathrm{Re}\iint \big(a_j(x)-a_j(y)\big) \frac{x_j-y_j}{
|x-y |^6}\big|u(t,x)\big|^2 \big|u(t,y) \big|^2dxdy
\\
& \approx 8\int_{|x|\leq R} |\nabla u(t,x)|^2 dx -8 \iint_{\Omega_1} \frac{|u(t,x)|^2|u(t,y)|^2}{|x-y|^4}dxdy \\
&\quad + O\big( \int_{|x|\approx R} \frac{|u(t,x)|^2}{R^2} dx +
\int_{|x|\approx R} |\nabla u(t,x)|^2 dx +  \iint_{\Omega_2}
\frac{|u(t,x)|^2|u(t,y)|^2}{|x-y|^4}dxdy  \big)\\
& \geq C_{ \delta_0} \int |\nabla u_0|^2dx,
\endaligned
\end{equation*}
where
\begin{equation*}
\aligned \Omega_1 & =\big\{ (x,y)\in \mathbb{R}^d\times
\mathbb{R}^d; |x|\leq R, |y|\leq R \big\};\\
\Omega_2&=\big\{(x,y)\in \mathbb{R}^d\times \mathbb{R}^d;
|x|\thicksim R \big\} \cup \big\{(x,y)\in \mathbb{R}^d\times
\mathbb{R}^d; |y|\thicksim R \big\}.\\
\endaligned
\end{equation*}

From (\ref{con6}) , we have
\begin{equation*}
\aligned C_{ \delta_0}  t  \int |\nabla u_0|^2dx \leq \big|
z'_R(t)-z'_R(0)\big| \leq2C  R^2 \int |\nabla u_0|^2dx.
\endaligned
\end{equation*}
We have a contradiction for $t$ large unless $u_0 \equiv 0$.

{\bf Proof of Theorem \ref{rigthm}: } It is analogue to the proof of
\cite{KeM06}, \cite{Mer01}. Assume that $u_0 \not\equiv 0$, then
\begin{equation}\label{posquan}
\aligned \int |\nabla u_0|^2 dx > 0.
\endaligned
\end{equation}
From Lemma \ref{usefulfacts}, we have
\begin{equation*}
\aligned E(u_0)\geq C_{\delta_0} \int |\nabla u_0|^2 dx > 0.
\endaligned
\end{equation*}
Because of Proposition \ref{strengthm}, we only need to consider the
case where there exists $\{t_n \}^{+\infty}_{n=1}$, $t_n \geq 0$, such
that
\begin{equation*}
\aligned \lambda(t_n) \rightarrow 0.
\endaligned
\end{equation*}
We claim that
\begin{equation*}
\aligned t_n \rightarrow T_{+}(u_0).
\endaligned
\end{equation*}
Indeed, if $t_n \rightarrow t_0 \in [0, T_{+}(u_0))$, then we have
for all $R>0$
\begin{equation*}
\aligned \int_{|x|>R}|v(t_n, x)|^{2^*}dx & = \int _{|x|>R}
\big|\frac{1}{\lambda(t_n)^{\frac{d-2}{2}}} u(t_n,
\frac{x}{\lambda(t_n)})\big|^{2^*}dx \\
& =  \int _{|x|>\frac{R}{\lambda(t_n)}} | u(t_n, x)|^{2^*}dx.
\endaligned
\end{equation*}
Because of $u\in C^0_t([0, T_{+}(u_0)); \dot{H}^1)$, we have
\begin{equation*}
\aligned \int_{|x|>R}|v(t_0, x)|^{2^*}dx & =0, \quad \forall\ R>0.
\endaligned
\end{equation*}
It is in contradiction with the fact that
\begin{equation*}
\aligned \int |\nabla v(t_0, x)|^2 dx = \int |\nabla u(t_0, x)|^2dx
>0.
\endaligned
\end{equation*}

Now after possibly redefining $\{t_n \}^{+\infty}_{n=1}$, we can
assume that
\begin{equation}\label{fact}
\aligned \lambda(t_n) \leq 2 \inf_{t\in [0, t_n]} \lambda(t).
\endaligned
\end{equation}
From the hypothesis, we have
\begin{equation*}
\aligned w_n(x)=\frac{1}{\lambda(t_n)^{\frac{d-2}{2}}}u\big(t_n,
\frac{x}{\lambda(t_n)}\big) \rightarrow w_0 \quad \text{in}\ \
\dot{H}^1.
\endaligned
\end{equation*}
By Proposition \ref{dlbv}, we have
\begin{equation*}
\aligned \int \big| \nabla w_n(x) \big|^2 dx &= \int \big| \nabla
u(t_n,x) \big|^2 dx < (1-\overline{\delta}) \int \big| \nabla W(x)
\big|^2 dx, \\
E(w_n) & =E(u(t_n))=E(u_0)<E(W).
\endaligned
\end{equation*}
Hence, we obtain
\begin{equation*}
\aligned \int \big| \nabla w_0 \big|^2 dx  & \leq
(1-\overline{\delta})
\int \big| \nabla W(x) \big|^2 dx \\
0< E(w_0)& =E(u_0) < E(W).
\endaligned
\end{equation*}
Thus $w_0 \not \equiv 0$. Let us now consider solutions
$w_n(\tau,x), w_0(\tau, x)$ of (\ref{equ1})  with data $w_n(x),
w_0(x)$ at $\tau=0$, defined in maximal intervals $\tau \in
(-T_{-}(w_n), 0]$ and $\tau \in (-T_{-}(w_0),0]$, respectively.
\begin{figure}[hb]
\qquad \qquad \qquad \qquad \quad \includegraphics[width=0.45\textwidth, angle=0]{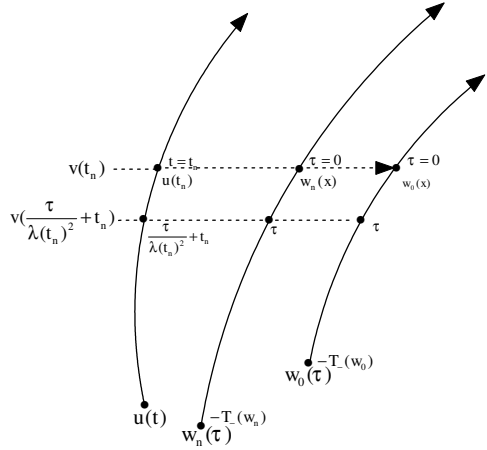}
\caption[]{A description of the normalization on $\lambda(t)$.}
\end{figure}

Since $w_n(x) \rightarrow w_0(x)$ in $\dot{H}^1$, we have from
Remark \ref{ct} that
\begin{equation}\label{prop}
\aligned \lim_{ \overline{n \rightarrow +\infty}} T_{-}(w_n) & \geq T_{-}(w_0),
\\
w_n(\tau,x) \rightarrow w_0(\tau,x) \quad & \text{in}\ \ \dot{H}^1,
\ \forall \ \tau \in (-T_{-}(w_0),0].
\endaligned
\end{equation}
By the uniqueness of solution of (\ref{equ1}), we have
\begin{equation*}
\aligned w_n(\tau, x)=
\frac{1}{\lambda(t_n)^{\frac{d-2}{2}}}u\big(\frac{\tau}{\lambda(t_n)^2}+t_n,
\frac{x}{\lambda(t_n)}\big), \quad \text{for}\ \
\frac{\tau}{\lambda(t_n)^2}+t_n \geq 0.
\endaligned
\end{equation*}
Now we claim that
\begin{equation}\label{taulimit}
\aligned \lim_{ \overline{n \rightarrow +\infty}} t_n \lambda(t_n)^2 \geq
T_{-}(w_0).
\endaligned
\end{equation}
Indeed, if not, then $\displaystyle \lim_{ \overline{n \rightarrow +\infty}} t_n
\lambda(t_n)^2 \rightarrow \tau_0 < T_{-}(w_0)$, from (\ref{prop}),
we have as $n\rightarrow +\infty$
\begin{equation*}
\aligned w_n(-t_n \lambda(t_n)^2,
x)=\frac{1}{\lambda(t_n)^{\frac{d-2}{2}}}u_0\big(\frac{x}{\lambda(t_n)}
\big) \rightarrow w_0(-\tau_0, x) \quad \text{in} \ \ \dot{H}^1.
\endaligned
\end{equation*}
Note that from $\lambda(t_n)\rightarrow 0$, we have as $n\rightarrow +\infty$
\begin{equation*}
\aligned
\frac{1}{\lambda(t_n)^{\frac{d-2}{2}}}u_0\big(\frac{x}{\lambda(t_n)}
\big) \rightharpoonup 0 \quad \text{in} \ \ \dot{H}^1,
\endaligned
\end{equation*}
thus we obtain that $w_0(-\tau_0)\equiv 0$, which yields a
contradiction.

From (\ref{taulimit}), we have that for fixed $\tau \in
(-T_{-}(w_0), 0]$ and sufficiently large $n$,
\begin{equation*}
\aligned 0 \leq \frac{\tau}{\lambda(t_n)^2}+t_n \leq t_n,
\endaligned
\end{equation*}
$v(\frac{\tau}{\lambda(t_n)^2}+t_n,x)$,
$\lambda(\frac{\tau}{\lambda(t_n)^2}+t_n)$ are defined and we have
\begin{equation*}
\aligned v(\frac{\tau}{\lambda(t_n)^2}+t_n,x) & =
\frac{1}{\lambda(\frac{\tau}{\lambda(t_n)^2}+t_n)^{\frac{d-2}{2}}}
u\big(\frac{\tau}{\lambda(t_n)^2}+t_n,
\frac{x}{\lambda(\frac{\tau}{\lambda(t_n)^2}+t_n)} \big) \\
& =
\frac{1}{\widetilde{\lambda}_n(\tau)^{\frac{d-2}{2}}} w_n\big(\tau,
\frac{x}{\widetilde{\lambda}_n(\tau)} \big),
\endaligned
\end{equation*}
where
\begin{equation*}
\aligned \widetilde{\lambda}_n(\tau)
=\frac{\lambda(\frac{\tau}{\lambda(t_n)^2}+t_n)}{\lambda(t_n)} \geq
\frac12
\endaligned
\end{equation*}
because of the fact (\ref{fact}). After passing to a subsequence, we
can assume that
\begin{equation*}
\aligned \widetilde{\lambda}_n(\tau) \rightarrow
\widetilde{\lambda}_0(\tau) \in [\frac12, +\infty].
\endaligned
\end{equation*}
Hence, we have
\begin{equation*}
\aligned
 v(\frac{\tau}{\lambda(t_n)^2}+t_n,x) \rightarrow \frac{1}{\widetilde{\lambda}_0(\tau)^{\frac{d-2}{2}}}w_0\big(\tau, \frac{x}{\widetilde{\lambda}_0(\tau)}
 \big)=v_0(\tau,x) \in \overline{K}.
\endaligned
\end{equation*}

Now we claim that
\begin{equation*}
\aligned \widetilde{\lambda}_0(\tau) < +\infty.
\endaligned
\end{equation*}
If not, from
\begin{equation*}
\aligned \frac{1}{\widetilde{\lambda}_n(\tau)^{\frac{d-2}{2}}}
w_n\big(\tau, \frac{x}{\widetilde{\lambda}_n(\tau)} \big)
\rightarrow \frac{1}{\lambda_0(\tau)^{\frac{d-2}{2}}}w_0\big(\tau,
\frac{x}{\lambda_0(\tau)} \big)= v_0(\tau,x),
\endaligned
\end{equation*}
we have
\begin{equation*}
\aligned w_0(\tau)=0,
\endaligned
\end{equation*}
which yields a contradiction.

So far, $w_0(\tau)$, $v_0(\tau)$ and $\widetilde{\lambda}_0(\tau)$
satisfy the conditions of Proposition \ref{strengthm}, we obtain
that
\begin{equation*}
\aligned w_0 \equiv 0,
\endaligned
\end{equation*}
which yields a contradiction. This completes the proof.

\begin{center}

\end{center}
\end{document}